\documentclass[10pt]{article}
\usepackage{amssymb,amsmath,amsfonts}
\usepackage{dsfont}
\usepackage{amsthm}
\usepackage{graphicx,color,subfig}
\usepackage{xspace}

\usepackage[latin1]{inputenc}
\usepackage[T1]{fontenc}

\usepackage{mathrsfs}
\usepackage{amscd}

\let \Im \relax
\DeclareMathOperator{\Im}{Im}
\newcommand{\Con}{\ensuremath{\mathscr C}}

\renewcommand{\S}{\ensuremath{\mathscr S}}

\newcommand{\est}[1]{\langle #1 \rangle}

\newcommand{\mb}[1]{\ensuremath{\mathbb{#1}}}
\newcommand{\N}{{\mb{N}}}

\newcommand{\R}{{\mb{R}}}
\newcommand{\C}{{\mb{C}}}

\newcommand{\eps}{\varepsilon}

\renewcommand{\d}{\ensuremath{\partial}}

\newenvironment{prooff}[1][Proof]{\textbf{#1.} }{\hfill \rule{0.5em}{0.5em}}

\DeclareMathOperator{\op}{op}
\DeclareMathOperator{\Op}{Op}

\DeclareMathOperator{\Ker}{ker}

\newcommand{\D}{\ensuremath{\mathscr D}}

\newcommand{\nhd}{neighborhood\xspace}

\usepackage[active]{srcltx}

\newcommand{\tno}[1]{ |\hskip -0,5pt |\hskip -0,5pt | #1 |\hskip -0,5pt |\hskip -0,5pt |}

\usepackage[english]{babel}

\newtheorem{lemma}{Lemma}[section]
\newtheorem{theorem}{Theorem}
\newtheorem{proposition}[lemma]{Proposition}

\theoremstyle{definition}

\newtheorem{remark}{Remark}

\makeatletter
\def\keywords{
    \vspace{1ex}
    \noindent
    \if@twocolumn
      \small{\bf  Keywords}\/---$\!$    \else
      \begin{center}\small\ {\bf Keywords}\end{center}\quotation\small
    \fi}
\def\endkeywords{\vspace{0.6em}\par\if@twocolumn\else\endquotation\fi
    \normalsize\rm}
\makeatother

\begin{document}

\title{Spectral analysis on interior transmission eigenvalues}

\author{ Luc Robbiano\thanks{Laboratoire de Math\'ematiques de Versailles, Universit\'e de Versailles St Quentin,
CNRS UMR 8100, 45, Avenue des \'Etats-Unis, 78035 Versailles, France. e-mail : Luc.Robbiano@uvsq.fr
}}

\maketitle
\begin{abstract}
In this paper we prove some results on interior transmission eigenvalues. First, under reasonable assumptions, we prove that the spectrum is a discrete countable set and the generalized eigenfunctions spanned a dense space in the range of resolvent.
This is a consequence of spectral theory of Hilbert-Schmidt operators. The main ingredient is to prove a smoothing property of resolvent. This allows to prove that a power of the resolvent is Hilbert-Schmidt.  
We obtain an estimate of the number of eigenvalues, counting with  multiplicities, with modulus less than $t^2$ when $t$ is large.
We prove also some estimate on the resolvent near the real axe when the square of the index of refraction 
is not real. 
Under some assumptions we obtain lower bound on the resolvent using the results obtained by Dencker, Sjöstrand and Zworski on the pseudospectra.

\end{abstract}


\tableofcontents

\section{Introduction}
In this paper we prove existence of infinite number of interior transmission eigenvalues under some condition on the index of refraction. Remind the problem. Let $\Omega$ a smooth bounded domain in $\R^n$. Let $n(x)$ a smooth function defined in $\overline \Omega$, called the index of refraction.
The problem is to find $k$ and $ (w,v)$ such that 
\begin{equation}\label{eq : ite}
 \left\{
	\begin{array}{l}
	 	\Delta w +k^2n(x) w=0 \text{ in } \Omega,\\ 
		\Delta v +k^2 v=0 \text{ in } \Omega,\\
		w=v \text{ on } \d \Omega,\\
		\d_\nu w=\d_\nu v  \text{ on } \d \Omega,
	\end{array}\right.
\end{equation}
where $\d_\nu$ is the exterior normal derivative to $\d\Omega$.
We consider here the function $n(x)$ complex valued. In the physical models, we have $n(x)=n_1(x)+in_2(x)/k$ where $n_j$ are real valued. Taking $u=w-v$ and $\tilde v=k^2v$, we obtain the following equivalent system
\begin{equation}\label{eq : ite Syl}
\left\lbrace 
 	\begin{array}{ll}
 		&\big( \Delta +k^2(1+m)\big) u+mv=0\text{ in } \Omega,\\
		&(\Delta +k^2)v=0 \text{ in } \Omega,\\
		&u=\partial_\nu u=0 \text{ on } \partial \Omega,
	\end{array}\right. 
\end{equation}
where, for simplicity, we have replaced $\tilde v$ by $v$ and $n$ by $1+m$.

Under some assumptions on $n(x)$, for instance our results apply for $n(x)=n_1(x)+in_2(x)/k$ where $n_j$ are real valued, we prove that the associated resolvent is compact on $H^2(\Omega)\oplus L^2(\Omega)$ (see Theorem~\ref{Resolv comp}) and we have a countable set of $k^2_j$ and generalized finite dimensional eigenspace $E_j$ such that $\cup_{j\in\N} E_j$ spanned a dense space in the range of the resolvent (see Theorem~\ref{Spect infini}). When $n(x) $ is real, Päivärinta ans Sylvester~\cite{PS08} have proved that there exist interior transmission eigenvalues, 
 Cakoni, Gintides and Haddar~\cite{CGH10} proved that the set of $k^2_j$ is infinite and discrete. For $n(x)$ complex valued  Sylvester~\cite{Sy12} proved that this set is discrete finite or infinite.
The case where $n(x)=1$ in a part of $\Omega$ that is the presence of cavities in $\Omega$ was considered by Cakoni, Çayören and Colton~\cite{CCC08},  Cakoni, Colton and Haddar~\cite{CCH10}. Here maybe because we use pseudo-differential calculus we do not have problems with cavities.
 In~\cite{HKOP10, HKOP11} Hitrik, Krupchyk, Ola and Päivärinta studied same type of problems where the Laplacian is replaced by an elliptic operator with constant coefficients of order $m\ge 2$. They proved in some cases, existence of interior transmission eigenvalue and the generalized eigenfunctions span a dense space. The proof uses the property of trace class operators and requires that $m>n$. Here as we consider the power of resolvent we have not this restriction and we consider the Laplacian for seek of simplicity but we can replace the Laplacian by a general elliptic operator of order 2 with real coefficients.

We can complete this result giving a weak version of Weyl law. If we denote by $N(t)$ the number of $|k_j|$, counting with multiplicities, smaller than $t$, we prove that $N(t)\le Ct^{n+4}$ (see Theorem~\ref{Fct compt}). 
In \cite{LV12}, \cite{LV11-arxiv} and \cite{LV12-arxiv},  Lakshtanov and Vainberg  study a problem as the \eqref{eq : ite} where the boundary condition $ \d_\nu w=\d_\nu v $ is replaced by $ \d_\nu w=a(x)\d_\nu v   $ where $a(x)\not=1$ for all $x\in\d\Omega$. In this case the problem is elliptic in the sense that if $(f,g)\in L^2(\Omega)$ then $(w,v)\in H^2\Omega)$ (see \eqref{syst. resolv. LV} for the definition of the system with force terms). In \cite{LV12bis-arxiv} they prove an lower bound on the counting function $N(t)$, if $n(x) $ is real and for the problem~\eqref{eq : ite}. To be precise they prove $N(t)\ge C t^n$ where $C>0$. Actually they consider only the real eigenvalues but it is not clear that the bound is sharp even for real eigenvalues.

In \cite{CK92} Colton and Kress prove that if $\Im n(x)\ge 0$ for all $x\in\overline\Omega$ and $\Im n\not=0$ then $k^2$ is not real. Here we give an estimate on the resolvent for $k^2\in\R$ (see Theorem~\ref{Maj Res}). This result is based on the Carleman estimates and following the same way as in the context of control theory, stabilization, scattering (see for instance~\cite{LR95}, \cite{LR97}, \cite{CR12}).

In some case (see Theorem~\ref{min Res}) we can give lower bound on the resolvent using the result obtained by 
 Dencker,  Sj{\"o}strand  and Zworski~\cite{DSZ04} on the pseudospectra. Even if the bounds obtained by Carleman's method and by the pseudospectra results have the same size, we cannot apply  both methods in the same situation.

The interest of the problem~\eqref{eq : ite} is related with the Theorem~8.9 in Colton and Kress~\cite{CK92} first proved by Colton, Kirsch and Päivärinta~\cite{CKP89}.
Here we give a quick survey of this result. Let $n(x)$ defined in $\Omega$ as the one in~\eqref{eq : ite} and by 1 in $\R^n\setminus\Omega$. We assume $n(x)\not=1$ in $\Omega$. For $k\in\R$, let $u$ the solution of the following problem
\begin{equation*}
 \left\{
	\begin{array}{l}
	 	\Delta u +k^2n(x) u=0 \text{ in } \R^n,\\ 
		u=u^i+u^s,\\
		\lim_{r\to +\infty}r^{(n-1)/2}(\frac{\d u^s}{\d r}-iku^s)=0,
	\end{array}\right.
\end{equation*}
where $u^i$ is a solution to $\Delta u^i+k^2u^i=0$ in $\R^n$, and $r=|x|$. We have $u^s(x)=\frac{e^{ik|x|}}{|x|^{(n-1)/2}}u_\infty(\hat x)+O\left( \frac{1}{|x|^{(n+1)/2}}\right) $, where $\hat x=x/|x|$.

If $u^i(x)=e^{ikxd}$, where $d\in\mb{S}^{n-1}$  we denote the corresponding $u_\infty$ by $u_\infty(\hat x,d)$. If we assume that $\Omega$ is connected and contains 0, the space spanned by $u_\infty(\hat x,d)$ for $d\in\mb{S}^{n-1}$ is dense in $L^2(\mb{S}^{n-1})$ if and only if the space of solution of~\eqref{eq : ite} is reduced to 0. This problem can be interpreted as follow, $u^i$ is a incident plane wave and $u_\infty$ is the first relevant term of $u$ created by the perturbation $n$ localized in $\Omega$.
This kind of result is interesting in the field of inverse  scattering problem. I let the reader interested by this field to find more information in~\cite{CK92} and in the recent survey given by Cakoni and Haddar~\cite{CH12}.

	\subsection{Results}

Let $\Omega $ a $\Con^\infty$ bounded domain in $\R^n$. Let $n(x)\in \Con^\infty(\overline\Omega)$ complex valued. We denote by $m(x)=n(x)-1$. We consider also the case where $n(x)=n_1(x)+in_2(x)/k$ where $n_j$ are real valued. We assume that for all $x\in\overline \Omega$, $n(x)\not=0$, or $n_1(x)\not=0$ or equivalently $m(x)\not= -1$. We assume there exists a \nhd $W$ of $\d \Omega$ such that for $x\in\overline W$,  $n(x)\not=1$ or  $n_1(x)\not=1$ or equivalently $m(x)\not= 0$. Actually if $n(x)\not=1$ for all $x\in\d\Omega$, such a \nhd $W$ exists.

Let $C_e$ the cone in $\C$ defined by 
\begin{equation}\label{def C}
 C_e=\{ z\in\C,\ \exists x\in\overline \Omega,\ \exists \lambda\ge 0,\text{ such that }z=-\lambda(1+\overline m(x))\}.
\end{equation}

In the case where $n=n_1+in_2/k$, $C_e=(-\infty,0]$ if $n_1(x)>-1$ for all $x\in\overline \Omega$, and $C_e=[0,+\infty)$ if $n_1(x)<-1$ for all $x\in\overline \Omega$.

Our regularity result will be given in the Sobolev spaces. We use the following notations.
 The $L^2(\Omega)$-norm  will be denoted by $\| \cdot \|$. For $s\in\R$,  we denote the usual $H^s$ norm in $\R^n$  by $\| w\|_{H^s}^2=\int (1+|\xi|^2)^s|\hat u(\xi)|^2d\xi$ where $\hat u$ is the classical Fourier transform. The $H^s$ space on $\Omega$ will be denoted by $\overline H^s(\Omega)$ and we say that $w$ a distribution in $\Omega$ is in $\overline H^s(\Omega)$ if there exists $\beta\in H^s$ such that $\beta_{|\Omega}=w$. The norm is given by $\| w\|_{\overline H^s(\Omega)}=\inf\{\| \beta\|_{H^s},\text{ where }\beta_{|\Omega}=w\}$. For $q\in\N$ it is classical that $\|w\|_{\overline H^q(\Omega)}^2$ is equivalent to $\sum_{|\alpha|\le q}\| \d^\alpha_x w\|^2$ (see for instance \cite[Vol. 3, Corollary B.2.5]{Hobo83}).

We denote by $B_z(u,v)=(f,g)$ the mapping defined from $H^2_0(\Omega)\oplus \{v\in L^2(\Omega), \ \Delta v\in L^2(\Omega)\}$ to $L^2(\Omega)\oplus L^2(\Omega)$  by

\begin{equation}\label{eq : resolvent 1}
\left\lbrace 
 	\begin{array}{ll}
 		&\big( \frac1{1+m}\Delta -z\big) u+\frac{m}{1+m}v=f\text{ in } \Omega\\
		&(\Delta -z)v=g \text{ in } \Omega\\
	\end{array}\right. 
\end{equation}

\begin{theorem}\label{th : Existence}
 Assume $C_e\not=\C$, then there exists $z\in \C$ such that $B_{z}$ is  bijective from $H^2_0(\Omega)\oplus \{v\in L^2(\Omega), \ \Delta v\in L^2(\Omega)\}$ to $ L^2(\Omega)\oplus  L^2(\Omega)$.
\end{theorem}

If for $z\in \C$ the solution $B_z(u,v)=(f,g)$ exists we  denote by $R_z(f,g)=(u,v)$.

\begin{theorem}\label{Resolv comp}
 Assume $C_e\not=\C$, there exists $z\in\C$ such that the resolvent $R_{z}$ from $\overline H^2(\Omega)\oplus L^2(\Omega)$ to itself is compact.

In particular, we can apply the Riesz theory, the spectrum is finite or a discrete countable set. If $\lambda\not=0$ is in the spectrum, $\lambda $ is a eigenvalue with finite  dimensional associated generalized eigenspace. 
\end{theorem}

\begin{remark}
 This result improves the Sylvester's theorem~\cite[theorem 2]{Sy12} with respect the geometrical assumption on $m$. Nevertheless, here the regularity  assumption  on $m$ is stronger than the assumption given in \cite{Sy12}.
\end{remark}

\begin{remark}
 Actually if $z_0 \not\in C_e\cup (-\infty,0]$ for all $\lambda $ large enough we can take $z=\lambda z_0$ in the Theorems~\ref{th : Existence} and \ref{Resolv comp}.
\end{remark}

\begin{remark}
 Actually we can also consider $R_z$ on $ L^2(\Omega)\oplus  L^2(\Omega)$ to itself. The range is in $H^2_0(\Omega)\oplus L^2(\Omega)$. Then with our regularity results we can prove that $R_z^2$ is a mapping from  $ L^2(\Omega)\oplus  L^2(\Omega)$ to $\overline{H}^4(\Omega)\oplus \overline{H}^2(\Omega)$. In particular $R_z^2$ is compact from $ L^2(\Omega)\oplus  L^2(\Omega)$ to itself and we can deduce the same properties on the spectrum of $R_z$ as in Theorem~\ref{Resolv comp}.
\end{remark}

In general for a non self-adjoint problem,  we cannot claim that the spectrum is non empty. In the following theorem, with a stronger assumption on $C_e$, we can prove that the spectrum is non empty.

We say that $C_e$ is contained in a sector with angle less than $\theta$ if there exist $\theta_1<\theta_2$, such that $C_e\subset \{ z\in \C, \ z=0 \text{ or } \frac{z}{|z|}=e^{i\varphi},\text{ where }\theta_1\le \varphi\le\theta_2\}$, and $\theta_2-\theta_1\le\theta $.

\begin{theorem}\label{Spect infini}
 Assume that $C_e$ is contained in a sector with angle less than $\theta$ with $\theta<{2\pi}/{p}$ and $\theta<\pi/2$ where $4p>n$. Then there exists $z$ such that the spectrum of $R_{z}$ is infinite and the space spanned by the generalized eigenspaces is dense in $H^2_0(\Omega)\oplus\{ v\in L^2(\Omega),\ \Delta v\in L^2(\Omega)\}$. 
 \end{theorem}
\begin{remark}
 This result is based on the theory given in Agmon~\cite{Ag65} and using the spectral results on Hilbert-Schmidt operators. In this theory we deduce that the spectrum is infinite from the proof that the generalized eigenspaces is dense in the closure of the range of $R_{z}$. Here we prove that $R_{z}^p$ is a Hilbert-Schmidt operator if $4p>n$. We can deduce the spectral decomposition of $R_{z}$ from the one of $R_{z}^p$.
\end{remark}

We can prove a weak Weyl law.
Let $z_j$ the elements of the spectrum of $R_{z}$ and $E_j$ the generalized associated eigenspace. We denote by $N(t)=\sum_{|z_j|^{-1}\le t^2}\dim E_j$.

\begin{theorem}\label{Fct compt}
 Under the same assumption as in Theorem~\ref{Spect infini}, then there exists $C>0$ such that $N(t)\le Ct^{4+n}$.
\end{theorem}

\begin{remark}
 I do not know if this result is optimal. The estimate is lower than the usual Weyl law which is in $t^{n}$. This is due to the estimates obtained on the resolvent which are  different than the one used to prove the usual Weyl law.
\end{remark}

Here we give some ideas to obtain the Theorems~\ref{Spect infini} and \ref{Fct compt} using the method given in Agmon~\cite{Ag65}.
First we prove a regularity result, that is we consider the iterate of $R_{z}$, we have $R_{z}^k$ is bounded from $H^2\oplus L^2$ to $H^{2k+2}\oplus H^{2k}$. This implies that $R_{z}^k$ is an Hilbert-Schmidt operator if $k$ is large enough and we can use the spectral theory for this operator class. The main problem to prove the regularity result is at the boundary. To do this we reduce the problem to the boundary by using the pseudo-differential calculus. It is well-known that for an elliptic problem, we can find a relation between the two traces of a solution. These two relations, for $u$ and $v$ in~\eqref{eq : resolvent 1} and the assumption on the two traces of $u$ allows to compute the trace of $v$ by the data. Actually the coupling is very weak because it involves  a lower order term, consequently  we obtain a weak estimate.

In the context of stabilization or control  for wave equation, there are a lot of results on decreasing of energy obtained by Carleman estimates (see for instance~\cite{LR95, LR97, CR12}). This method allows to give quantitative results related with uniqueness result. We use here the same method to prove an estimate on the resolvent near the real axe for a complex index of refraction.
The  theorem below is an quantitative version of the Theorem~8.12 given by Colton and Kress~\cite{CK92}. Here it is more convenient to use the variables introduce in \eqref{eq : ite}. Let $(w,v)$ solutions of

\begin{equation}\label{syst. resolv. LV}
 \left\{
	\begin{array}{l}
	 	\Delta w +k^2n(x) w=f \text{ in } \Omega\\ 
		\Delta v +k^2 v=g \text{ in } \Omega\\
		w=v \text{ on } \d \Omega\\
		\d_\nu w=\d_\nu v  \text{ on } \d \Omega,
	\end{array}\right.
\end{equation}
where $(f,g)\in L^2(\Omega)\oplus L^2(\Omega)$. We denote $\tilde R_{k^2}(f,g)=(w,v)$. Remark that $\tilde R_{k^2}$ exists except for a discrete set of $k^2$. Indeed we can check that $R_{-k^2}(f-g,k^2g)=(w-v,k^2v)$ this gives $(w,v)$ if $(-k^2)^{-1}$ is not in the spectrum of $R_0$.
\begin{theorem}\label{Maj Res}
 We assume that $\Im n\ge 0$ and $\Im n\not\equiv 0$ or if $n(x)=n_1(x)+in_2(x)/k$, $n_2(x)\ge 0$ and $n_2\not\equiv0$. Then there exist $C_1>0$ and $C_2>0$ such that
$\| \tilde R_{k^2}\|\le C_1e^{C_2|k|}$ for all $k\in\R$. Here $\| \cdot\| $ denote the norm of the operator from $ L^2(\Omega)\oplus L^2(\Omega)$ to itself.
\end{theorem}

\begin{remark}
 We have the same result if we assume $\Im n\le 0$ and $\Im n\not\equiv 0$ or if $n(x)=n_1(x)+in_2(x)/k$, $n_2(x)\le 0$ and $n_2\not\equiv0$.
\end{remark}

In the context of non self-adjoint operator, the spectrum is not the most relevant notion. Actually Davies~\cite{Da99} introduced the notion of the pseudospectrum. Roughly speaking this set is defined by the points $z$ where the resolvent is large. This notion is related with the notion of ill-conditioned for the matrix. Here we use the result proved by Dencker, Sjöstrand and Zworski~\cite{DSZ04} to obtain a lower bound on the norm of resolvent.

\begin{theorem}\label{min Res}
 Assume there exist $\xi_0\in\R^n$ and $x_0\in \Omega$ such that $\Im((\overline{n}(x_0))\est{\xi_0,\d_xn(x_0)})\not=0$. Then for all $N>0$, $ \sup \{ |k|^{-N}\| \tilde R_{k^2}\|, k^2=r(n(x_0))^{-1}, r>0\}=+\infty$. Moreover if $n$ is an analytic function in a \nhd of $x_0$ there exists $C>0$ such that  $ \sup \{ e^{-C|k|}\|\tilde R_{k^2}\|, k^2=r(n(x_0))^{-1}, r>0\}=+\infty$.
\end{theorem}

\begin{remark}
 Even if the lower bound, in analytic case, is of the same type of the upper bound obtained in Theorem~\ref{Maj Res}, we cannot apply  both theorems for the same direction $z$. Actually $k^2$ is in general in $\C$, if we want apply the Theorem~\ref{Maj Res} we need $n(x_0)\in\R$, that is $\Im n(x_0)=0$, as $\Im n(x)\ge0$ then we cannot have $\Im \d_x n(x_0)\not=0$ to apply the Theorem~\ref{min Res}.
\end{remark}

\subsection{Ouline}

In the Section~\ref{sec reg} we prove the main technical results. Roughly speaking if the data $(f,g)$ are more regular in $H^s$ norm we prove that the solution $(u,v)$ is also more regular. More precisely, we prove for $p\ge 0$, that if $(f,g)\in \overline H^{2p+2}(\Omega)\oplus  \overline H^{2p}(\Omega)$ then $R_{z}(f,g)=(u,v)\in \overline H^{2p+4}(\Omega)\oplus  \overline H^{2p+2}(\Omega)$. This proves first that $R_{z}$ is compact as an operator from $ \overline H^{2}(\Omega)\oplus  L^2(\Omega)$ to itself and $R^p_{z}$ is an operator from $ \overline H^{2}(\Omega)\oplus  L^2(\Omega)$ to $\overline H^{2p+2}(\Omega)\oplus  \overline H^{2p}(\Omega)$. This implies that $R^p_{z}$ is an Hilbert-Schmidt operator on $ \overline H^{2}(\Omega)\oplus  L^2(\Omega)$ if $4p>n$. To prove the regularity result, first we prove an estimate on $u$ in the subsection~\ref{sub u}. It is an easy estimate to obtain as $u_{|\d\Omega}=0$, $u$ satisfies a classical Dirichlet problem. Second we prove in Subsection~\ref{sub v int} the regularity of $v$ in all compact in $\Omega$. As $v$ satisfies an elliptic equation, far away the boundary of $\Omega$ it is a classical result. In the third Subsection~\ref{sub v bord} we prove the regularity result on $v$ in a \nhd of $\d\Omega$. The idea to do this is to explain $v$ by the unknown traces of $v$. This description allows to obtain a relation between $v_{|\d\Omega}$ and $\d_n v_{|\d\Omega}$. Then we can use this formula on $v$ in the equation on $u$. The fact that $u_{|\d\Omega}=0$ and $\d_n u_{|\d\Omega}=0$  gives  another relation between $v_{|\d\Omega}$ and $\d_n v_{|\d\Omega}$. These relations allow to determine $v_{|\d\Omega}$ and $\d_n v_{|\d\Omega}$ with $(f,g)$. This explicit formula, in sense of pseudo-differential calculus, allows to prove the regularity result. We need also following the same way to prove an estimate of the $L^2$ norm of $v$ by the $L^2$ norm of $f$. This implies a weak convergence result. Actually the problem is that $v$ has the same regularity than $f$. In particular if we consider the resolvent as a operator from $L^2(\Omega)\oplus L^2(\Omega)$ to itself, we cannot prove that the resolvent is compact. Here we avoid this problem by the assumption that $f\in\overline H^2(\Omega)$.

In the Section~\ref{sec spect} we recall some result proved in Agmon~\cite{Ag65}  and we apply this to prove the Theorems~\ref{Resolv comp}, \ref{Spect infini} and \ref{Fct compt}.  

In the Section~\ref{sec est res} we prove some a priori bound on the resolvent. In Subsection~\ref{sub upper} we prove an upper bound on the resolvent near the real axe when the imaginary part of the refraction index have a sign and is not identically null. 
The main tool is to use the  interpolation estimate obtained from the Carleman estimate. In the Subsection~\ref{sub lower} we use the result  obtained by Dencker, Sjöstrand and Zworski~\cite{DSZ04} on the pseudospectra to obtain lower bound on the resolvent. Roughly speaking this result says that when the operator is not elliptic in the semi-classical sense, even if a point is not in the spectrum, the resolvent is generically large. 

As we use deeply the semi-classical pseudo-differential calculus in the Section~\ref{sec reg},
in the Appendix~\ref{App pseudo}  we fix the notation used in the rest of paper, we recall the classical results used, we give some ideas of proof on the action of the pseudo-differential operators  on $\overline H^s(\Omega)$ spaces, we give  some computations on integrals to obtain some explicit  formulae used in Subsection~\ref{sub v bord}. This allows to give the explicit first term of the resolvent in sense of the semi-classical pseudo-differential calculus.

\section{Regularity results}\label{sec reg}

We describe now the idea of the proof.
As we want prove an estimate when $|z|$ is large we will compute the resolvent in semi-classical framework. We multiply the Equations \eqref{eq : resolvent 1} by $h^2$, we denote by $\mu=- h^2z$ where $\mu$ belongs to a bounded domain of $\C$, $a=1/(1+m)$ and $V=m/(1+m)$. We change $(f,g)$ in $(-f,-g)$.

Remind the assumption on $m$, we have $m(x)\not=-1$ for all $x\in\overline \Omega$ and $m(x)\not=0$ for $x$ in a \nhd of $\d\Omega$.

Thus  following \eqref{eq :  resolvent 1}, we obtain the system
\begin{equation}\label{eq : resolvent sc}
\left\lbrace 
 	\begin{array}{ll}
 		&\big(-ah^2 \Delta -\mu\big) u-h^2Vv=h^2f\text{ in } \Omega\\
		&(-h^2\Delta -\mu)v=h^2g \text{ in } \Omega\\
		&u=\partial_\nu u=0 \text{ on } \partial \Omega.
	\end{array}\right. 
\end{equation}
The goal of this section is to prove the following estimates if $s\ge 0$. The result is given using the semi-classical $H^s$ norm, see Appendix~\ref{App pseudo} for the definition of these spaces.

\begin{theorem}\label{th est spect}
We assume that 
 for all $x\in\overline \Omega$, and all $\xi\in\R^n$, $a(x)|\xi|^2-\mu\not=0$ and
 $|\xi|^2-\mu\not=0$. Let
  $s\ge 0$, there exists $h_0>0$ such that for 
$f\in\overline H^{2+s}_{sc}(\Omega)$,
 $g\in \overline H^s_{sc}(\Omega)$, $u\in \overline H^{1+s}_{sc}(\Omega)\cap H^{2}_0(\Omega)$ and $v\in \overline H^s_{sc}(\Omega)$ solutions of system \eqref{eq : resolvent sc}
then $u\in \overline  H^{4+s}_{sc}(\Omega)$ , $v\in \overline H^{2+s}_{sc}(W)$ and for $h\in (0,h_0)$ we have,

\begin{align}
 &\| u\|_{\overline{H}^{4+s}_{sc}(\Omega)}\lesssim  h^2\| f\|_{\overline{H}^{2+s}_{sc}(\Omega)}+h^4 \| g\|_{\overline{H}^s_{sc}(\Omega)} \label{norme H s+2 sur u}
 \\
& \| v \|_{\overline{H}^{2+s}_{sc}(\Omega)}\lesssim  \| f\|_{\overline{H}^{2+s}_{sc}(\Omega)}+h^2 \| g\|_{\overline{H}^s_{sc}(\Omega)}.\label{norme H s+1 sur v}
\end{align}
\end{theorem}

First we prove an estimate on $u$. For this, we work globally in $\Omega$. The estimate on $v$ is more difficult to obtain. In a first step we prove an estimate in the interior by usual pseudo-differential tools. In a second step we prove the estimate in a \nhd of the boundary $\d \Omega$ and we finish the proof. This will be do in the following three sections.

In the proof we use semi-classical pseudo-differential calculus. We give in Appendix~\ref{App pseudo}, the results used, trace formula, action of pseudo-differential operators on Sobolev space,  the parametrices.

\subsection{Estimate on $u$}\label{sub u}

The goal of this section is to proof a weak version of \eqref{norme H s+2 sur u}.
\begin{lemma}\label{lem est u}
We assume that for all $x\in\overline \Omega$, and all $\xi\in\R^n$, $a(x)|\xi|^2-\mu\not=0$. There exists $h_0>0$ such that for $s\ge 0$, for all
 $f\in\overline H^{s}_{sc}(\Omega)$, $v\in \overline H^{s}_{sc}(\Omega)$, and $u\in H^{2}_0(\Omega)$ solution of 
\begin{equation*}
\left\lbrace 
 	\begin{array}{ll}
 		&\big(-ah^2 \Delta -\mu\big) u=h^2Vv+h^2f\text{ in } \Omega\\
		&u=\partial_\nu u=0 \text{ on } \partial \Omega,
	\end{array}\right. 
\end{equation*}
 then $u\in \overline  H^{2+s}_{sc}(\Omega)$ and for $h\in(0,h_0)$,
\begin{equation}\label{est lem u}
 \| u\|_{\overline H^{2+s}_{sc}(\Omega)}\lesssim
h^2\| f\|_{\overline H^{s}_{sc}(\Omega)}+h^2\|v \|_{\overline H^{s}_{sc}(\Omega)}.
\end{equation}
\end{lemma}

\begin{prooff}
 As $u\in H^2_0(\Omega)$, we can extend $u$ by 0 in the exterior of $\Omega$ and $u$ satisfied the same equation in the whole space. Here we extend also $v$ and $f$ by 0, this makes sense at least in $L^2$. We have 
\begin{equation}\label{eq sur u}
 \big(-ah^2 \Delta -\mu\big)\underline u=h^2V\underline v+h^2 \underline f\text{ in } \R^n,
\end{equation}
where we denote, for $w\in L^2(\Omega)$, by 
$$ \underline w(x)=\left\lbrace 
\begin{array}{ll}
 &w(x) \textrm{ if } x\in \Omega\\
& 0 \textrm{ if } x\not\in \Omega.
\end{array}\right. 
$$

Let $N\ge s+2$, we take a parametrix $Q$ of $-ah^2\Delta-\mu$, this is possible because we assume  for all $x\in \overline\Omega$, $a(x)|\xi|^2-\mu\not=0$. We have $Q(-ah^2\Delta-\mu)=\chi+hK$ where $\chi\in\Con^\infty_0(\R^n) $ , $\chi=1$ in a \nhd of $\Omega$, $K$ is of order $-N$
and $Q$ is of order $-2$.

Applying $Q$ to Equation~\eqref{eq sur u} we obtain,
\begin{equation*}
 \underline u+hK \underline u=h^2Q \underline f+h^2 Q\left( V\underline v\right) 
\end{equation*}
As $Q$ is a mapping on  the Sobolev spaces (see \eqref{eq ope para}) we obtain,
\begin{equation*}
 \| u\|_{\overline H^{2+s}_{sc}(\Omega)}\lesssim h\|u\|_{L^2(\Omega)}+h^2\| f\|_{\overline H^{s}_{sc}(\Omega)}+h^2\|v \|_{\overline H^{s}_{sc}(\Omega)}.
\end{equation*}
We can absorb the term $h\|u\|_{L^2(\Omega)} $ by the left hand side and this imply \eqref{est lem u}.  
\end{prooff}

\subsection{Estimate on $v$ in interior of $\Omega$}\label{sub v int}

To prove an estimate on $v$ in interior of $\Omega$ we follow essentially the same way than  in the proof of the estimate on $u$ given  in the  previous section except that we cannot extend $v$ in exterior of $\Omega$ but we use semi-classical pseudo-differential calculus in open relatively compact in $\Omega$.
The estimate proved is given in the following lemma.
\begin{lemma}\label{lem int v}
We assume  for  all $\xi\in \R^n$, we have $|\xi|^2-\mu\not=0$.
Let $\chi\in\Con^\infty_0(\R^n)$ supported in $W$ relatively compact in $\Omega$, and $s\ge 0$, there exists $h_0>0$ such that for 
 $g\in \overline H^s_{sc}(\Omega)$ and $v\in \overline H^{s+1}_{sc}(\Omega)$ solution of 
\begin{equation*}
 (-h^2\Delta -\mu)v=h^2g \text{ in } \Omega,
\end{equation*}
then $v\in \overline H^{2+s}_{sc}(W)$ and for $h\in (0,h_0)$ we have,
\begin{equation}\label{est v int}
 \| \chi v\|_{ H^{2+s}_{sc}}\lesssim h\| v\|_{\overline H^{s+1}_{sc}(\Omega)}+h^2\| g  \|_{\overline H^{s}_{sc}(\Omega)}.
\end{equation}
\end{lemma}

\begin{prooff}
 In the sequel we can take $\chi$ supported in $W$ or $\chi=1$ on $W$ and supported in a compact of $\Omega$. We can essentially repeat the proof of Lemma~\ref{lem est u} given to estimate $u$.

As we have assumed that  $|\xi|^2-\mu\not=0$, we can  take a parametrix $\tilde Q$ of $(-h^2\Delta-\mu)$ defined globally in $\R^n$, such that we have $\tilde Q(-h^2\Delta-\mu)=\chi+hK$ where $K$ is of order $-1$,
$\tilde Q$ is of order $-2$. Let $\chi_j\in\Con_0^\infty(\R^n)$,   $j=1,2$
supported in a compact of $\Omega$ where $\chi_1=1$ on the support of $\chi$ and $\chi_2 =1$ on support of $\chi_1$. By pseudo-differential calculus we have $\tilde Q\chi_1(-h^2\Delta-\mu)=\chi +hK_{-1}$ where  $K_{-1} $ is of order $-1$.
Now we have $\tilde Q\chi_1(-h^2\Delta-\mu)\chi_2=\chi +hK_{-1}\chi_2$. We can localized the equation on $v$ and we have  $\chi_1(-h^2\Delta-\mu)\chi_2 v=h^2\chi_1 g$.
Applying $\tilde Q$ to this equation we obtain
\begin{equation}\label{Est v}
\chi v+hK_{-1}( \chi_2 v)=h^2\tilde Q( \chi_1 g).
\end{equation}
Taking the $H^{2+s}_{sc}$ norm of $\chi v$ we obtain \eqref{est v int}.
\end{prooff}

\subsection{Estimate on $v$ in a \nhd of the boundary }\label{sub v bord}

\begin{prooff}[Proof of Theorem~\ref{th est spect}]

Taking account the Lemma~\ref{lem est u} and \ref{lem int v}, to acheive the proof,  we need an estimate on $v$ near the boundary $\d \Omega$.
It is  well-known that we can find in a \nhd $W$ of the boundary $\partial \Omega$ a system of coordinates such that the Laplacian can be written $\partial_{x_n}^2+R(x,\partial_{x'})+\alpha(x)\partial_{x_n}$, where $x'$ are the coordinates on the manifold $\partial \Omega$, $x_n\in(0,\eps)$, $\Psi(W)=\partial \Omega\times (0,\eps)$, $\Psi$ is the change of coordinates and $R$ is a differential operator on $\partial \Omega$ of order 2 depending of the parameter $x_n$.

We keep the notation $u$, $v$, $a$ in the coordinates $x$ instead of $u\circ \Psi$, etc. The Equations \eqref{eq : resolvent sc} become
\begin{equation}\label{eq : resolvent scl}
\left\lbrace 
 	\begin{array}{ll}
 		&\big(a(D_{x_n}^2+R(x,D')+h\alpha D_{x_n}) -\mu\big) u-h^2Vv=h^2f\text{ in } 	\partial D\times (0,\eps)\\
		&(D_{x_n}^2+R(x,D')+h\alpha D_{x_n} -\mu\big)v=h^2g \text{ in } \partial D\times (0,\eps)\\
		&u=\partial_\nu u=0 \text{ on } \partial \Omega\times \{ 0\}.
	\end{array}\right. 
\end{equation}
We have taken the semi-classical notations, $D_{x_j}= \frac{h}{i}\partial_{x_j}$, $D'=(D_{x_1},\cdots, D_{x_{n-1}})$. Actually $R$ and $V$ may depend explicitly of $h$ but this introduces no problem in the estimates, if $V_{|h=0}\not=0$, this is the case if $n=n_1+in_2/k$, we have $V_{|h=0}(x)=(1-n_1)/n_1\not=0$ by assumption.

Let $w\in L^2(x_n>0)$, we denote by
\begin{equation}\label{prolongement 0}
 	\underline w=1_{x_n>0}w=\left\lbrace 
		\begin{array}{l}
 			w \text{ if } x_n>0\\
			0 \text{ if } x_n<0.
		\end{array}\right. 
\end{equation}
Usually we use $\underline w$ but some time it is more convenient to use $ 1_{x_n>0}w $. We have 
\begin{align*}
& D_{x_n}\underline w=1_{x_n>0}D_{x_n}w+\frac{h}{i}w_{|x_n=0}\otimes\delta_{x_n=0}\\
&D_{x_n}^2\underline w=1_{x_n>0}D_{x_n}^2w+\frac{h}{i}D_{x_n}w_{|x_n=0}\otimes \delta_{x_n=0}+\frac{h}{i}w_{|x_n=0}\otimes D_{x_n}\delta_{x_n=0}.
\end{align*}
Here and in the sequel, for simplicity we denote by $w_{|x_n=0}$ the limit, when $x_n$ goes to 0 with $x_n>0$, of $w(x',x_n)$, if the limit exists. Here the distributions $u$ and $v$ are solutions of elliptic equations then the limits exist in a space of distributions.

From \eqref{eq : resolvent scl}, we obtain
\begin{equation}\label{eq : resolvent bord}
\left\lbrace 
 	\!\!\!\!\!\!\!\!\begin{array}{ll}
 		&\big(a(D_{x_n}^2+R(x,D')+h\alpha D_{x_n}) -\mu\big) \underline u-h^2V\underline v=h^2\underline f\text{ in } 	\partial \Omega\times (-\eps,\eps)\\
		&(D_{x_n}^2+R(x,D')+h\alpha D_{x_n} -\mu\big)\underline v=h^2\underline g +\frac{h}{i}\gamma_0\otimes\delta_{x_n=0}+\frac{h}{i}\gamma_1\otimes D_{x_n}\delta_{x_n=0}\text{ in } \partial \Omega\times (-\eps,\eps)  ,
	\end{array}\right. 
\end{equation}
where $\gamma_0=D_{x_n}v_{|x_n=0}+\alpha h v_{|x_n=0}$ and $ \gamma_1=v_{|x_n=0}$. We can consider these equations for $x_n \in (-\eps,\eps)$ indeed the coefficients of $R$ are smooth up the boundary, and we can extend $R$  in a \nhd of the boundary for $x_n<0$. The functions $ \underline u$ and $ \underline v$ are null for $x_n<0$ so the equations are relevant only to take account the boundary terms. Remark,  in the first equation because the traces of $u$ are null, they are not boundary terms.

The main goal of this section, is to obtain estimates on $\gamma_0$ and $\gamma_1$.

Now we search, using the equations \eqref{eq : resolvent bord}, two relations between the traces of $v$. 
First we localize $v$ in a \nhd of the boundary. We denote by $w=\chi_0 v$ where $\chi_0\in \Con^\infty(\R)$, 
$\chi_0(x_n)=1$ in a \nhd of boundary, for instance if $|x_n|\le \eps/4$ and $\chi_0(x_n)=0$ if $|x_n|\ge \eps/2$. 
From the second equation of \eqref{eq : resolvent bord} we obtain
\begin{equation}\label{eq : sur w}
 (D_{x_n}^2+R(x,D')+h\alpha D_{x_n} -\mu\big)\underline w=h^2\chi_0\underline g +\frac{h}{i}\gamma_0\otimes
 \delta_{x_n=0}+\frac{h}{i}\gamma_1\otimes D_{x_n}\delta_{x_n=0} +hK\underline v\text{ on } \partial \Omega\times \R ,
\end{equation}
where $K$ is a first order differential operator coming from the commutator between $D_{x_n}^2 $ or $D_{x_n}$ and $\chi_0$.

Let $\chi_1\in\Con^\infty(\R)$ such that $\chi_1\chi_0=\chi_0$ for instance $\chi_1(x_n)=1$  if $|x_n|\le \eps/2$ and $\chi_0(x_n)=0$ if $|x_n|\ge 3\eps/4$. By assumption we have $\xi_n^2+R(x,\xi')-\mu\not= 0$ then by semi-classical pseudo-differential calculus there exists $\tilde Q$ of order -2 such that $\tilde Q(D_{x_n}^2+R(x,D')+h\alpha D_{x_n} -\mu\big)=\chi_1 +h\tilde K$ where $\tilde K$ is of order $-N$ where $N\ge s+2$.
Applying $\tilde Q$ to \eqref{eq : sur w}, we obtain
\begin{align}
 &\chi_1 \underline{w}=\underline{w}=\tilde Q(\frac{h}{i}\gamma_0\otimes\delta_{x_n=0}+\frac{h}{i}\gamma_1\otimes D_{x_n}\delta_{x_n=0}) +g_1\notag\\
&\text{where }g_1=-h\tilde K\underline w+h^2\tilde Q(\chi_0 \underline g)+h\tilde QK\underline v \text{ thus} \notag\\
&\| g_1\|_{\overline H^{s+2}_{sc}(\Omega)}\lesssim h\| v\|_{\overline H^{s+1}_{sc}(\Omega)}+h^2\| g\|_{\overline H^{s}_{sc}(\Omega)},
 \label{eq sur w 2}
\end{align}
actually we can estimate $\tilde K\underline w$ because $\underline w\in L^2(\R^n)$ and $\tilde K$ is smoothing. By this trick we have not to verify that $\tilde K$ is a mapping on  the $\overline H^s_{sc}$. In appendix~\ref{App pseudo}, Estimate~\eqref{eq ope para}, we have proved that a parametrix as $\tilde Q$ is a mapping on  the $\overline H^s_{sc}(\Omega)$.

By ellipticity assumption on $\xi_n^2+R(x,\xi')-\mu$, there exist $\rho_1(x,\xi')$ and $\rho_2(x,\xi')$ with $\Im \rho_1>0$ and $\Im\rho_2<0$ such that $\xi_n^2+R(x,\xi')-\mu=(\xi_n-\rho_1(x,\xi')) (\xi_n-\rho_2(x,\xi'))$ (see \ref{sb : roots param.}).

From~\eqref{rest term bord} and  Lemma~\ref{lem trace}  we have

\begin{align}
&\tilde Q(\frac h i\gamma\otimes D_{x_n}^k\delta_{x_n=0})=\op (\tilde q)\gamma\\
& \text{where } \tilde q=\frac 1 {2i\pi}\int_\R e^{ix_n\xi_n/h}\frac{\xi_n^k}{(\xi_n-\rho_1)(\xi_n-\rho_2)}d\xi_n
+r_{-2+k},
\end{align}
where we denote here and in the sequel  by $r_j$  an operator of order $j$. 

From Lemma~\ref{lem sym 1} we obtain if we restrict \eqref{eq  sur  w 2} on $\{ x_n=0\}$ and as $\underline w_{|x_n=0}=\gamma_1$, 

\begin{align}
\gamma_1=\op(\frac 1{\rho_1-\rho_2})\rho_0+\op(\frac{\rho_1}{\rho_1-\rho_2})\gamma_1+h\op(r_{-2})\gamma_0
+h\op(r_{-1})\gamma_1+(g_1)_{|x_n=0}.
\end{align}
Thus we obtain
\begin{align}\label{eq traces 0}
 &\op \left( \frac{ \rho_2}{\rho_2-\rho_1}\right) \gamma_1+\op \left( \frac{1}{\rho_2-\rho_1}\right) \gamma_0=(g_1)_{|x_n=0}+
h\op(r_{-2})\gamma_0+h\op(r_{-1})\gamma_1.
\end{align}

Then applying $ \op(\rho_2-\rho_1) $ on  both sides of \eqref{eq traces 0}, using estimate \eqref{eq sur w 2} to estimate $g_1$, trace formula and by pseudo-differential calculus we obtain
\begin{align}
 &\op(\rho_2)\gamma_1+\gamma_0=g_2  \textrm{ where } g_2=\op(\rho_2-\rho_1) (g_1)_{|x_n=0}+
h\op(r_{-1})\gamma_0+h\op(r_{0})\gamma_1
\notag\\
&\text{and } | g_2|_{H^{s+1/2}_{sc}(\d \Omega)}\lesssim h^{1/2}\| v\|_{\overline H^{s+1}_{sc}(\Omega)}+h^{3/2}\| g\|_{\overline H^{s}_{sc}(\Omega)}+h| \gamma_0|_{H^{s-1/2}_{sc}(\d \Omega)}+h| \gamma_1|_{H^{s+1/2}_{sc}(\d \Omega)}  \label{eq traces 1}
\end{align}

To obtain a second equation on the traces, we use the first equation of \eqref{eq : resolvent bord}. As before there exist $Q$ of order -2, $K$ of order $-N-4$, $\chi_2$ such that $\chi_2\chi_0=\chi_0$ and such that $Q(a(D^2_{x_n}+R+ h\alpha D_{x_n})-\mu)=\chi_2+hK$ .

We apply $Q$ to the first equation of \eqref{eq : resolvent bord}, we obtain
\begin{align}
 &\chi_2 \underline u=h^2Q (V\underline w)+g_3 \notag \\
&\text{where } g_3=-hK\underline u+h^2Q (V(1-\chi_0)\underline v)    +h^2Qf  \notag \\
&\| g_3\|_{\overline H^{s+4}_{sc}(\Omega)}\lesssim h\| u\|_{\overline H^{s+2}_{sc}(\Omega)}+h^2\|          (1-\chi_0)\underline v  \|_{\overline H^{s+2}_{sc}(\Omega)}
+h^2\| f\|_{\overline H^{s+2}_{sc}(\Omega)}.  \label{eq u}
\end{align}
Using Lemma~\ref{lem int v} we can estimate the term $(1-\chi_0)\underline v $, using the Lemma~\ref{lem est u} we can estimate $ \| u\|_{\overline H^{s+2}_{sc}(\Omega)} $ and we obtain,
\begin{equation*}
 \| g_3\|_{\overline H^{s+4}_{sc}(\Omega)}\lesssim h^3\|          v  \|_{\overline H^{s+1}_{sc}(\Omega)}+h^4\| g\|_{\overline H^{s}_{sc}(\Omega)}+h^2\| f\|_{\overline H^{s+2}_{sc}(\Omega)}.
\end{equation*}

We replace $\underline w$ by the Formula \eqref{eq sur w 2}, we take the trace on $x_n=0$, as $u_{|x_n=0}=0$, we obtain
\begin{align}
 &\left[ Q\left(V\left( \tilde Q(\frac{h}{i}\gamma_0\otimes\delta_{x_n=0}+\frac{h}{i}\gamma_1\otimes D \delta_{x_n=0})
\right) 
 \right) \right]_{|x_n=0}=g_4,\label{eq bord 2}\\
&\text{where }g_4=-h^{-2}(g_3)_{|x_n=0}-[QV(g_1)]_{|x_n=0},\text{ then }  \notag \\
&| g_4|_{\overline H^{s+7/2}_{sc}(\d \Omega)}\lesssim h^{1/2}\| v\|_{\overline H^{s+1}_{sc}(\Omega)}+h^{3/2}\| g\|_{\overline H^{s}_{sc}(\Omega)}+h^{-1/2}\| f\|_{\overline H^{s+2}_{sc}(\Omega)}.\label{est g4}
\end{align}
As $a(x)(\xi_n^2+R(x,\xi'))-\mu$ is elliptic, the polynomial in $\xi_n$ has two roots $\lambda_j$ one satisfies $\Im \lambda_1>0$  and the other $\Im \lambda_2<0$ (see section~\ref{sb : roots param.}).
 We have $a(x)(\xi_n^2+R(x,\xi'))-\mu=a(\xi_n-\lambda_1)(\xi_n-\lambda_2)$. The principal symbol of $Q$ is $\frac{\chi_2}{a(\xi_n-\lambda_1)(\xi_n-\lambda_2)}$. The principal symbol of $QV\tilde Q$ is $\frac{V\chi_1\chi_2}{a(\xi_n-\lambda_1)(\xi_n-\lambda_2)(\xi_n-\rho_1)(\xi_n-\rho_2)}$.

Following the same method used to obtain~\eqref{eq traces 0} from~\eqref{eq sur w 2}, we have by \eqref{rest term bord}, Lemmas~\ref{lem trace} and \ref{lem sym 2}, as $\chi_2=1$ in \nhd of $\d \Omega$
\begin{align*}
& \op\left( \frac{V(\lambda_2-\lambda_1+\rho_2-\rho_1)}{a(\lambda_1-\lambda_2)(\lambda_1-\rho_2)(\rho_1-\lambda_2)(\rho_1-\rho_2)}
\right) \gamma_0\notag \\
&\qquad \qquad +
\op\left( \frac{V(\rho_2\lambda_2-\rho_1\lambda_1)}{a(\lambda_1-\lambda_2)(\lambda_1-\rho_2)(\rho_1-\lambda_2)(\rho_1-\rho_2)}
\right) \gamma_1=g_5  \notag \\
&\text{where } g_5=g_4+h\op( r_{-4})\gamma_0+h\op( r_{-3})\gamma_1.
\end{align*}
As $V\not=0$ in a \nhd of $\d \Omega$, we can apply $\op(\frac a V (\lambda_1-\lambda_2)(\lambda_1-\rho_2)(\rho_1-\lambda_2)(\rho_1-\rho_2) )$, we obtain
\begin{align}
& \op\left((\lambda_2-\lambda_1+\rho_2-\rho_1)
 \right) \gamma_0
+ \op\left( {\left( \rho_2\lambda_2-\rho_1\lambda_1\right) }
\right) \gamma_1=g_6, \notag \\
&\text{where } g_6=\op(\frac a V (\lambda_1-\lambda_2)(\lambda_1-\rho_2)(\rho_1-\lambda_2)(\rho_1-\rho_2) )g_5+h\op(r_{0})\gamma_0 +h\op(r_{1})\gamma_1   \notag\\
&\text{then } | g_6|_{H^{s-1/2}_{sc}(\d \Omega)}\lesssim
h^{1/2}\| v\|_{\overline H^{s+1}_{sc}(\Omega)}+h^{3/2}\| g\|_{\overline H^{s}_{sc}(\Omega)}+h^{-1/2}\| f\|_{\overline H^{s+2}_{sc}(\Omega)}\notag\\
&\hskip 3,4cm+h| \gamma_0 |_{H^{s-1/2}_{sc}(\d \Omega)} +h| \gamma_1 |_{H^{s+1/2}_{sc}(\d \Omega)}\label{eq trace 2}
\end{align}
Now we have two equations on traces  by \eqref{eq traces 1} we can replace $\gamma_0$ in \eqref{eq trace 2}. We obtain by pseudo-differential calculus
\begin{align}
 \op(\rho_2\lambda_2-\lambda _1\rho_1)\gamma_1-\op((\lambda_2-\lambda_1+\rho_2-\rho_1)\rho_2)\gamma_1=g_6-\op(\lambda_2-\lambda_1+\rho_2-\rho_1)g_2=g_7.\notag
\end{align}
This implies,
\begin{align}
 &\op((\rho_2-\rho_1)(\lambda_1-\rho_2))\gamma_1=g_7, \text{ with}\notag\\
&| g_7 |_{H^{s-1/2}_{sc}(\d \Omega)}\lesssim 
h^{1/2}\| v\|_{\overline H^{s+1}_{sc}(\Omega)}+h^{3/2}\| g\|_{\overline H^{s}_{sc}(\Omega)}+h^{-1/2}\| f\|_{\overline H^{s+2}_{sc}(\Omega)}\notag\\
&\hskip 2,5cm +h| \gamma_0 |_{H^{s-1/2}_{sc}(\d \Omega)} +h| \gamma_1 |_{H^{s+1/2}_{sc}(\d \Omega)}. \label{est g7}
\end{align}
As $\Im \lambda_1>0$, $\Im \rho_1>0$ and $\Im \rho_2<0$, the symbol $(\rho_2-\rho_1)(\lambda_1-\rho_2)$ is elliptic, by inversion we obtain
\begin{align}
| \gamma_1 |_{H^{s+3/2}_{sc}(\d \Omega)}&\lesssim  h^{1/2}\| v\|_{\overline H^{s+1}_{sc}(\Omega)}+
h^{3/2}\| g\|_{\overline H^{s}_{sc}(\Omega)}+h^{-1/2}\| f\|_{\overline H^{s+2}_{sc}(\Omega)}\notag\\
&\quad+h| \gamma_0 |_{H^{s-1/2}_{sc}(\d \Omega)} +h| \gamma_1 |_{H^{s+1/2}_{sc}(\d \Omega)}.\label{est gamma 1}
\end{align}
and using \eqref{eq traces 1}, we obtain
\begin{align}
| \gamma_0 |_{H^{s+1/2}_{sc}(\d \Omega)}&\lesssim 
h^{1/2}\| v\|_{\overline H^{s+1}_{sc}(\Omega)}+h^{3/2}\| g\|_{\overline H^{s}_{sc}(\Omega)}+h^{-1/2}\| f\|_{\overline H^{s+2}_{sc}(\Omega)}\notag\\
&\quad+h| \gamma_0 |_{H^{s-1/2}_{sc}(\d \Omega)} +h| \gamma_1 |_{H^{s+1/2}_{sc}(\d \Omega)}.\label{est gamma 0}
\end{align}
Summing  \eqref{est gamma 1} and \eqref{est gamma 0} we have for $h_0$ small enough
\begin{align}
 | \gamma_1 |_{H^{s+3/2}_{sc}(\d \Omega)}+| \gamma_0 |_{H^{s+1/2}_{sc}(\d \Omega)}\lesssim
h^{1/2}\| v\|_{\overline H^{s+1}_{sc}(\Omega)} +h^{3/2}\| g\|_{\overline H^{s}_{sc}(\Omega)}
+h^{-1/2}\| f\|_{\overline H^{s+2}_{sc}(\Omega)} .  \notag
\end{align}

From \eqref{eq  sur  w 2} and from estimate~\eqref{app est int par trace} obtained in Appendix~\ref{App pseudo}
\begin{align}
 \| w\|_{\overline H^{s+2}_{sc}(\Omega)}&\lesssim \| g_1\|_{\overline H^{s+2}_{sc}(\Omega)}+h^{1/2}\left( | \gamma_0|_{H^{s+1/2}_{sc}(\d \Omega)}+| \gamma_1|_{H^{s+3/2}_{sc}(\d \Omega)}\right) \notag\\
&\lesssim 
h\| v\|_{\overline H^{s+1}_{sc}(\Omega)}
+h^{2}\| g\|_{\overline H^{s}_{sc}(\Omega)}+\| f\|_{\overline H^{s+2}_{sc}(\Omega)}.\label{est w}
\end{align}
We can now estimate $v$, we have by \eqref{est w}, Lemma~\ref{lem int v},
\begin{align*}
 \| v\|_{\overline H^{s+2}_{sc}(\Omega)}&\le \| w\|_{\overline H^{s+2}_{sc}(\Omega)}+ \| (1-\chi_0 )v\|_{\overline H^{s+2}_{sc}(\Omega)}\notag\\
&\lesssim h\| v\|_{\overline H^{s+1}_{sc}(\Omega)}
+h^{2}\| g\|_{\overline H^{s}_{sc}(\Omega)}+\| f\|_{\overline H^{s+2}_{sc}(\Omega)}.
\end{align*}
This implies 
$$
\| v\|_{\overline H^{s+2}_{sc}(\Omega)}\lesssim h^{2}\| g\|_{\overline H^{s}_{sc}(\Omega)}+\| f\|_{\overline H^{s+2}_{sc}(\Omega)},
$$
if $h_0$ is small enough.
Using this estimate and \eqref{est lem u} we obtain
$$
\| u\|_{\overline H^{4+s}_{sc}(\Omega)}\lesssim
h^2\| f\|_{\overline H^{s+2}_{sc}(\Omega)}+h^4\| g\|_{\overline H^{s}_{sc}(\Omega)}.
$$
These two last estimates imply the Theorem~\ref{th est spect}.
\end{prooff}

\begin{remark}
Actually in this proof we have assumed $v\in \overline H^{s+1}_{sc}(\Omega)$. To prove the result with $v\in \overline H^{s}_{sc}(\Omega)$, we argue in two steps, first following the same proof we can obtain  $v\in \overline H^{s+1}_{sc}(\Omega)$ and second the proof given above gives $v\in \overline H^{s+2}_{sc}(\Omega)$
\end{remark}

The estimate proved above are not enough on $v$ with respect $f$. 
For the sequel we need to estimate $v$ by $f$ in $H^s$ norm.

\begin{proposition}\label{est v L2}

 We assume that 
 for all $x\in\overline \Omega$, and all $\xi\in\R^n$, $a(x)|\xi|^2-\mu\not=0$ and
 $|\xi|^2-\mu\not=0$. There exist $h_0>0$ and $C>0$ and for all $\delta>0$, let $\chi_\delta\in \Con^\infty$ supported in a $\delta$-\nhd of $\d \Omega$, there exists $C_\delta>0$
such that   for $h\in (0,h_0)$,
$f\in \overline{H}^s_{sc}(\Omega)$, $g=0$,
  $u\in \overline H^{s+1}_{sc}(\Omega)\cap H^2_0(\Omega)$, $v\in \overline{H}^s_{sc}(\Omega)$ and $\Delta v\in L^2(\Omega)$, solutions of system \eqref{eq : resolvent 1} we have the estimate,
\begin{align*} 
 \| v \|_{ \overline{H}^s_{sc}(\Omega)}\le C \| \chi_\delta f\|_{ \overline{H}^s_{sc}(\Omega)}+C_\delta h  \| f \|_{ \overline{H}^s_{sc}(\Omega)}. 
\end{align*}
\end{proposition}

\begin{prooff}
 The proof follows the previous one. We give only the modifications to do.
From \eqref{eq : sur w} we obtain \eqref{eq sur w 2} with the estimate 
\begin{equation}\label{est g1}
 \| g_1\|_{\overline H^{s+1}_{sc}(\Omega)}\lesssim h\| v \|_{\overline H^s_{sc}(\Omega)}.
\end{equation}
Thus we obtain \eqref{eq traces 0} and \eqref{eq traces 1} with the estimate 
\begin{align}
 | g_2|_{H^{s-3/2}(\d \Omega)}
&\lesssim | \op(r_1) (g_1)_{|\d\Omega}|_{H^{s-1/2}(\d \Omega)}+h|  \gamma_0|_{H^{s-5/2}(\d\Omega)} +h| \gamma_1|_{H^{s-3/2}(\d\Omega)}\notag\\
&\lesssim h^{1/2}\|  v \|_{\overline H^s_{sc}(\Omega)}+h|  \gamma_0|_{H^{s-5/2}(\d\Omega)} +h| \gamma_1|_{H^{s-3/2}(\d\Omega)}.\label{est g2}
\end{align}
We must modify \eqref{eq u} to obtain the term $\chi_\delta$.  We take the same $Q$ as in~\eqref{eq u} and we apply $Q\chi_\delta$ to the first equation from~\eqref{eq : resolvent bord}. We obtain 
\begin{align*}
Q\chi_\delta \big(a(D_{x_n}^2+R(x,D')+h\alpha D_{x_n}) -\mu\big) \underline u-h^2Q\chi_\delta V\underline v=h^2Q\chi_\delta \underline f\text{ in } 	\partial \Omega\times (-\eps,\eps).
\end{align*}
We have $\chi_\delta \big(a(D_{x_n}^2+R(x,D')+h\alpha D_{x_n}) -\mu\big)= \big(a(D_{x_n}^2+R(x,D')+h\alpha D_{x_n}) -\mu\big)\chi_\delta+L_1$ where $L_1$ is a differential operator of order 1 depending of $\delta$.
 As $Q(a(D^2_{x_n}+R+ h\alpha D_{x_n})-\mu)=\chi_2+hK_{-N}$  where $K_{-N}$ is of order $-N$, with $N\ge s+2$, and $ \chi_2\chi_\delta=\chi_\delta$, we have
\begin{align*}
 &\chi_\delta \underline u=h^2Q (\chi_\delta V\underline w)+g_3\\
&\text{where } g_3=
-hQL_1 \underline u+hK_{-N} \underline u+h^2Q(\chi_\delta f)\\
&\| g_3\|_{\overline H^{s+2}_{sc}(\Omega)}\le C_\delta h\| u\|_{\overline H^{s+1}_{sc}(\Omega)}
+Ch^2\|\chi_\delta f\|_{\overline H^{s}_{sc}(\Omega)} . 
\end{align*}
We have used that $Q$ and $L_1$, a differential operator, act on the $\overline H^s_{sc}$.
We can estimate $u$ by \eqref{est lem u}, this gives
\begin{equation*}
 \| g_3\|_{\overline H^{s+2}_{sc}(\Omega)}\le C_\delta h^3\| f\|_{\overline H^{s}_{sc}(\Omega)}+C_\delta h^3\| v \|_{\overline H^{s}_{sc}(\Omega)}
+Ch^2\|\chi_\delta f\|_{\overline H^{s}_{sc}(\Omega)} 
\end{equation*}
We replace $\underline w$ by its value given by the formula~\eqref{eq sur w 2} with the estimate~\eqref{est g1}. 
We obtain~\eqref{eq bord 2} with
\begin{equation}\label{est g4bis}
 | g_4|_{H^{s+3/2}_{sc}(\d\Omega)}  \le
C_\delta h^{1/2}\| f\|_{\overline H^{s}_{sc}(\Omega)}+C_\delta h^{1/2}\| v \|_{\overline H^{s}_{sc}(\Omega)}
+Ch^{-1/2}\|\chi_\delta f\|_{\overline H^{s}_{sc}(\Omega)}.
\end{equation}
If we compare this estimate with~\eqref{est g4} we see that the bad power of $h$ in front of $f$ is only 
a part of $f$ localized in a \nhd of the boundary.

We can follow the proof and we obtain~\eqref{eq trace 2} where the estimate on $g_6$  is
\begin{align*}
 | g_6|_{H^{s-5/2}(\d\Omega)}&\le C_\delta h^{1/2}\| f\|_{\overline H^{s}_{sc}(\Omega)}+C_\delta h^{1/2}\| v \|_{\overline H^{s}_{sc}(\Omega)}
+Ch^{-1/2}\|\chi_\delta f\|_{\overline H^{s}_{sc}(\Omega)}\\
&\quad +Ch|\gamma_0|_{H^{s-5/2}(\d\Omega)}+Ch|\gamma_1|_{H^{s-3/2}(\d\Omega)}.
\end{align*}
We have the Formula~\eqref{est g7} where $g_7$, from \eqref{est g2} is estimated by
\begin{align*}
 | g_7|_{H^{s-5/2}(\d\Omega)}&\le C_\delta h^{1/2}\| f\|_{\overline H^{s}_{sc}(\Omega)}+C_\delta h^{1/2}\| v \|_{\overline H^{s}_{sc}(\Omega)}
+Ch^{-1/2}\|\chi_\delta f\|_{L^2(\Omega)}\\
&\quad +Ch|\gamma_0|_{H^{s-5/2}(\d\Omega)}+Ch|\gamma_1|_{H^{s-3/2}(\d\Omega)}.
\end{align*}
By ellipticity and Formula~\eqref{est g7} we obtain 
\begin{align*}
 | \gamma_1|_{H^{s-1/2}(\d\Omega)}&\le C_\delta h^{1/2}\| f\|_{\overline H^{s}_{sc}(\Omega)}+C_\delta h^{1/2}\| v \|_{\overline H^{s}_{sc}(\Omega)}
+Ch^{-1/2}\|\chi_\delta f\|_{\overline H^{s}_{sc}(\Omega)}\\
&\quad +Ch|\gamma_0|_{H^{s-5/2}(\d\Omega)}+Ch|\gamma_1|_{H^{s-3/2}(\d\Omega)},
\end{align*}
and by~\eqref{eq traces 1} where $g_2$ satisties~\eqref{est g2}, we have,
\begin{align*}
  | \gamma_0|_{H^{s-3/2}(\d\Omega)}&\le C_\delta h^{1/2}\| f\|_{\overline H^{s}_{sc}(\Omega)}+C_\delta h^{1/2}\| v \|_{\overline H^{s}_{sc}(\Omega)}
+Ch^{-1/2}\|\chi_\delta f\|_{\overline H^{s}_{sc}(\Omega)}\\
&\quad +Ch|\gamma_0|_{H^{s-5/2}(\d\Omega)}+Ch|\gamma_1|_{H^{s-3/2}(\d\Omega)}.
\end{align*}
Summing the previous estimates and for $h$ small enough, we obtain,
\begin{align*}
 | \gamma_1|_{H^{s-1/2}(\d\Omega)}+ | \gamma_0|_{H^{s-3/2}(\d\Omega)}&\le C_\delta h^{1/2}\| f\|_{\overline H^{s}_{sc}(\Omega)}+C_\delta h^{1/2}\| v \|_{\overline H^{s}_{sc}(\Omega)}
+Ch^{-1/2}\|\chi_\delta f\|_{\overline H^{s}_{sc}(\Omega)}.
\end{align*}
From~\eqref{eq sur w 2} with $g_1$ satisfying~\eqref{est g1}, we have by \eqref{app est int par trace}
\begin{align}
 \| w\|_{\overline H^{s}_{sc}(\Omega)}&\le C \| g_1\|_{\overline H^{s}_{sc}(\Omega)}+Ch^{1/2}\left( | \gamma_0|_{H^{s-3/2}_{sc}(\d \Omega)}+| \gamma_1|_{H^{s-1/2}_{sc}(\d \Omega)}\right) \notag\\
&\le
C_\delta h\| v\|_{\overline H^{s}_{sc}(\Omega)}
+C_\delta h
\| f\|_{\overline H^{s}_{sc}(\Omega)}
+C\|\chi_\delta f\|_{\overline H^{s}_{sc}(\Omega)}.
\label{est wbis}
\end{align}
Using the formula \eqref{Est v} in the proof of Lemma~\ref{lem int v} with  $g=0$, we obtain $ \| \chi v\|_{\overline H^s_{sc}(\Omega)}\le C \|  v\|_{\overline H^s_{sc}(\Omega)}$.
This estimate and \eqref{est wbis} give
\begin{equation*}
  \| v\|_{\overline H^{s}_{sc}(\Omega)}\le
C_\delta h\| v\|_{\overline H^{s}_{sc}(\Omega)}
+C_\delta h
\| f\|_{L^2(\Omega)}
+C\|\chi_\delta f\|_{\overline H^{s}_{sc}(\Omega)}.
\end{equation*}
This estimate is also true for a fixed $\delta$ then we have for $h\in(0,h_0)$, $h_0$ small enough
\begin{equation*}
 \| v\|_{\overline H^{s}_{sc}(\Omega)}\le C \| f\|_{\overline H^{s}_{sc}(\Omega)}.
\end{equation*}
This with the previous estimate implies the Proposition~\ref{est v L2}.
\end{prooff}

\section{Existence and compactness}

In this section we prove the Theorems~\ref{th : Existence} and  \ref{Resolv comp}. 

\subsection{Proof of Theorem~\ref{th : Existence}}

\begin{prooff}
 We follow the proof given by Sylvester~\cite[Proposition 10]{Sy12}, we prove that le range of $B_{z}$ is closed and dense. 

To prove the range is closed we apply the a priori estimates prove in section~\ref{sec reg}. We recall that $a={1}/{(1+m)}$ and $V={m}/{(1+m)}$.

We remark that  if  we have $C_e\not= \C$ then $C_e\cup (-\infty,0] \not=\C$. Indeed, as $\overline \Omega$ is compact, $C_e$ is closed. If $C_e\cup (-\infty,0]=\C$ then $C_e\setminus  (-\infty,0]=\C\setminus (-\infty,0]$ as $\C\setminus  (-\infty,0]$ is dense in $\C$, we have $C_e=\C$.
 Let  $z_0$ such that $z_0\not\in C_e\cup  (-\infty,0]$,
  we can choose $|z_0|=1$, let $z=h^{-2}z_0$ we have $\mu=-z_0$. First we can estimate $\|v\|_{L^2(\Omega)}$ by Proposition~\ref{est v L2} with $s=0$ and $\delta$ fixed if $g=0$ and by Theorem~\ref{th est spect} if $f=0$. 
There exists $C>0$ such that  for all $|z|$ large enough, 
\begin{equation}\label{est sur v}
 \|v\|_{L^2(\Omega)}\le C\|f\|_{L^2(\Omega)}+\frac{C}{|z|^2}\|g\|_{L^2(\Omega)}.
\end{equation}
We can apply the Lemma~\ref{lem est u} with $s=0$, we obtain with the previous estimate on $v$,

\begin{equation}\label{est sur u}
 |z|^2\|u\|_{L^2(\Omega)}+|z|\| u\|_{\overline{H}^1(\Omega)}+\|u\|_{\overline{H}^2(\Omega)}\le C\|f\|_{L^2(\Omega)}+\frac{C}{|z|^2}\|g\|_{L^2(\Omega)}.
\end{equation}
Clearly these estimates prove that the range of $B_z$ is closed where the norm on the domain of the operator is given by the $H^2$ norm for $u$ and by $\|v\|+\|\Delta v\|$ for $v$.

To prove the density of the range of $B_z$ we prove that the orthogonal of the range is $\{ 0\}$. We recall the Green formula, if $v$ and $q$ are smooth functions in $\overline{\Omega}$ we have 
\begin{equation}\label{Green}
 (v|\Delta q)-(\Delta v|q)=(v|\partial_\nu q)_0-(\partial_\nu v|q)_0,
\end{equation}
where $(\cdot|\cdot)$ is the inner product in $\Omega$, $(\cdot|\cdot)_0$ is the inner product on $\partial \Omega$ and $\partial_\nu$ is the exterior normal derivative on $\partial \Omega$. Actually \eqref{Green} is true if $v$ smooth, $q\in L^2(\Omega)$  and  $\Delta q\in L^2(\Omega)$. Indeed, in this case it is well known that $q_{|\d \Omega}\in H^{-1/2}(\d\Omega)$ and $\d_\nu q_{|\d\Omega}\in H^{-3/2}(\d\Omega)$ then we can find $f_n$ and $g_n$ smooth functions such that $f_n$ goes to $\Delta q$ in $L^2(\Omega)$ and $g_n$ goes to $q_{|\d \Omega}$. Let $q_n$ the solution of $\Delta q_n=f_n$ in $\Omega$ and $(q_n)_{|\d\Omega}=g_n$, $q_n$ is a smooth function and by continuity $(\d_\nu q_n)_{|\d\Omega}$ goes to $\d_\nu q_{\d\Omega}$ in $H^{-3/2}(\d\Omega)$. Then we can pass to the limit in \eqref{Green}.

Let $p,q\in L^2(\Omega)$ and $u,v$ smooth functions in $\Omega$, if $(p,q)$ is in the orthogonal of the range we have
\begin{equation}\label{eq duale}
 (\Delta u-z(1+m)u+mv|p)+(\Delta v-zv|q)=0.
\end{equation}
We take $u,v\in\Con^\infty_0(\Omega)$ in \eqref{eq duale}, by integrating by part in distribution sense we have
\begin{align}
 \Delta p-\bar z(1+\bar m)p=0\text{ in }\Omega\label{eq sur p}\\
\Delta q-\bar zq+\bar m p=0\text{ in }\Omega.\label{eq sur q}
\end{align}
In particular $\Delta p$ and $\Delta q$ are in $L^2(\Omega)$, then we can apply \eqref{Green} to integrate by part in \eqref{eq duale} if now $u,v$ are smooth functions up the boundary with $u_{|\d\Omega}=\d_\nu u_{|\d\Omega}=0$. Using \eqref{eq sur p} and \eqref{eq sur q} we have 
\begin{equation*}
 (v|\d_\nu q)_0-(\d_\nu v|q)_0=0.
\end{equation*}
As $v_{|\d\Omega}$ and $\d_\nu v_{|\d\Omega}$ are arbitrary, we obtain $q_{|\d\Omega}=\d_\nu q_{|\d\Omega}=0$. By $\eqref{eq sur q}$, $q$ satisfies a Dirichlet boundary value problem and $\d_\nu q_{|\d\Omega}=0$, then $q\in\overline{H}^2_0(\Omega)$. By \eqref{eq sur p}, $\Delta p\in L^2(\Omega)$ and $p\in L^2(\Omega)$. We deduce that $(q,p)\in  H^2_0(\Omega)\oplus \{v\in L^2(\Omega), \ \Delta v\in L^2(\Omega\}$ satisfies the same kind of equation as $(u,v)$. Then the inequalities  \eqref{est sur v} and \eqref{est sur u} prove that $p=q=0$. This acheives the proof of Theorem~\ref{th : Existence}.
\end{prooff}

\subsection{Proof of Theorem~\ref{Resolv comp}}

\begin{prooff}
 We take the same $z$ as in the proof of Theorem~\ref{th : Existence}. We can apply Theorem~\ref{th est spect} with $s=0$ , we obtain in classical norm 
\begin{align}
 &|z|^2\|u\|_{L^2(\Omega)}+|z|\| u\|_{\overline{H}^1(\Omega)}+\|u\|_{\overline{H}^2(\Omega)}+
\frac1{|z|}\|u\|_{\overline{H}^3(\Omega)}+
\frac1{|z|^2}\|u\|_{\overline{H}^4(\Omega)}
\le C\|f\|_{\overline{H}^2(\Omega)}+\frac{C}{|z|^2}\|g\|_{L^2(\Omega)},\notag\\
&\|v\|_{L^2(\Omega)}+\frac{C}{|z|}\|v\|_{\overline{H}^1(\Omega)}+
\frac1{|z|^2}\|v\|_{\overline{H}^2(\Omega)}
\le C\|f\|_{\overline{H}^2(\Omega)}+\frac{C}{|z|^2}\|g\|_{L^2(\Omega)}.\label{Est norm R}
\end{align}
This proves that $R_z : \overline{H}^2(\Omega)\oplus L^2(\Omega)\to \overline{H}^4(\Omega)\oplus
 \overline{H}^2(\Omega)$,  then $R_z$ is compact from $\overline{H}^2(\Omega)\oplus L^2(\Omega)$ to itself. We can apply the Riesz theory.
\end{prooff}

\section{Spectral results}\label{sec spect}

Here we prove how the regularity results obtained in the  section~\ref{sec reg} allow to prove the spectral results. Actually the result obtained in Theorem~\ref{Resolv comp} is not enough to prove that the spectrum is a countable set. The theory given in Agmon~\cite{Ag65}  is based on the spectral decomposition of Hilbert-Schmidt operators.
We adapt two results given in Agmon to our case, the Lemma~\ref{lem est H-S} and the Proposition~\ref{prop infinie vp}. Following these results we will prove the Theorems~\ref{Spect infini} and~\ref{Fct compt}.

Let $T$ an  Hilbert-Schmidt operator  from $\overline H^2(\Omega)\oplus L^2(\Omega)$ to itself. We denote by $\tno{T}$ the Hilbert-Schmidt norm.  Let $(\phi_j)_{j\in\N}$ a Hilbert basis on $\overline H^2(\Omega)$ and $(\psi_k)_{k\in\N}$ a Hilbert basis on $L^2(\Omega)$, then $((\phi_j,0))_{j\in\N},((0,\psi_k))_{k\in\N}$ is a Hilbert basis on $\overline H^2(\Omega)\oplus L^2(\Omega)$. Let $T(\phi_j,0)=u_j=(u^0_j,u^1_j)$ and $T(0,\psi_k)=v_k=(v^0_k,v^1_k)$.
With these notations, we have  $\tno{T}^2= \sum_{j=0}^\infty ( \|u^0_j \|_{\overline H^2(\Omega)}^2+\|  v^0_j\|_{\overline H^2(\Omega)}^2
+\| u^1_j \|_{L^2(\Omega)}^2+\|  v^1_j\|_{L^2(\Omega)}^2)$.

We denote by $\|T\|_j$  the operator  norm from $\overline H^2(\Omega)\oplus L^2(\Omega)\to \overline H^{j+2}(\Omega)\oplus \overline H^{j}(\Omega)$, where $\overline H^{0}(\Omega)=L^2(\Omega)$.

\begin{lemma}\label{lem est H-S}
Let $m>n/2$, there exists $C>0$ such that if $T$ is a bounded operator from  $\overline H^2(\Omega)\oplus L^2(\Omega)\to\overline H^{m+2}(\Omega)\oplus \overline H^{m}(\Omega)$, then $T$ is a Hilbert-Schmidt operator and 
$$
\tno{T} \le C\|T\|^{n/(2m)}_m \|T\|^{1-n/(2m)}_0 .
$$
\end{lemma}

\begin{prooff}
 We follow the proof given  by Agmon~\cite[Theorem 13.5]{Ag65}.

 Let $u=T(\sum_{j=0}^Na_j(\phi_j,0))+T(\sum_{j=0}^Nb_j(0,\psi_j))=
\sum_{j=0}^Na_ju_j+\sum_{j=0}^Nb_jv_j=(u^0,u^1)$. We have  $u_j=(u^0_j,u^1_j)$ and $v_j=(v^0_j,v^1_j)$. 
We treat the term $u^0$.

If  $m>n/2$, $\overline H^m(\Omega)\subset L^\infty (\Omega)$ and (see~\cite[Lemma 13.2]{Ag65}), 
for $\alpha\in\N^n$, $|\alpha|\le 2$,  there exists $C>0$ such that  
$$
\|\d^\alpha v\|_{L^\infty(\Omega)} \le  C \| v\|_{\overline H^{m+|\alpha|}(\Omega)}^{n/(2m)}
\|v\|_{\overline H^{|\alpha|}(\Omega)}^{1-n/(2m)} 
$$
By the property on $T$, we have 
$$
\|u^0\|_{\overline H^{m+2}(\Omega)}^2\le C\|T\|^2_m\sum_{j=0}^N(|a_j|^2+|b_j|^2)\text{ and }
 \|u^0\|_{\overline H^{2}(\Omega)}^2\le C\|T\|^2_0\sum_{j=0}^N(|a_j|^2+|b_j|^2).
$$ 
Let $K=\|T\|^{n/(2m)}_m \|T\|^{1-n/(2m)}_0$ thus we have for $x\in\Omega$,
 and for $\alpha\in\N^n$, $|\alpha|\le 2$,  $|\d^\alpha u^0(x) |^2\le CK^2\sum_{j=0}^N(|a_j|^2+|b_j|^2)$. 
 We have $\d^\alpha u^0(x)=\sum_{j=0}^N(a_j\d^\alpha u^0_j(x)+b_j\d^\alpha v^0_j(x))$, 
 we  take in the previous inequality
 $a_j=\d^\alpha \bar u_j^0(x)$ and $b_j=
\d^\alpha\bar v_j^0(x)$,  
we sum on $\alpha$, we obtain for all $x\in\overline \Omega$,
\begin{align*}
 \sum_{|\alpha|\le 2}\left( \sum _{j=0}^N(|\d^\alpha u_j^0(x)|^2+|\d^\alpha v_j^0(x)|^2)\right) ^2 
\le CK^2\sum_{|\alpha|\le 2}\sum _{j=0}^N(|\d^\alpha u_j^0(x)|^2+|\d^\alpha v_j^0(x)|^2).
\end{align*}
 Thus 
$\sum_{|\alpha|\le 2}(|\d^\alpha u_j^0(x)|^2+|\d^\alpha v_j^0(x)|^2))\le CK^2$, integrating this on $\Omega$ 
(which is bounded) we find 
$\sum_{j=0}^N ( \|u^0_j \|_{\overline{H}^2(\Omega)}^2+\|  v^0_j\|_{\overline{H}^2(\Omega)}^2)\le CK^2$. 
As the right hand side does not depend on $N$ we can let $N$ goes to infinity.
We can treat by the same method the terms $\sum_{j=0}^N ( \|u^1_j \|_{L^2(\Omega)}^2+\|  v^1_j\|_{L^2(\Omega)}^2)$,
it suffices to repeat the previous proof without the derivative terms. This means that $\tno{T}$ 
is bounded by $CK=C\|T\|^{n/(2m)}_m \|T\|^{1-n/(2m)}_0$.
\end{prooff}

We give here a small improvement of the Theorem~16.4 in \cite{Ag65}. 

We introduce some notations. The inner product in $\overline{H}^2(\Omega)\oplus L^2(\Omega)$ will be denoted by $(\cdot|\cdot)$. Let $T$ an operator from $\overline{H}^2(\Omega)\oplus L^2(\Omega)$ to itself, if $\lambda^{-1}$  is in the resolvent set of $T$, we set $T_\lambda=T(I-\lambda T)^{-1}$. We remark that if $T$ is the resolvent of $P$, that is $PT=I$, then $T_\lambda$ is the resolvent of $P-\lambda I$. Indeed, 
$$
(P-\lambda I)T_\lambda=(P-\lambda I)T(I-\lambda T)^{-1}=(I-\lambda T)^{-1}-\lambda T(I-\lambda T)^{-1}=(I-\lambda T)^{-1}(I-\lambda T)=I.
$$
\begin{proposition}\label{prop infinie vp}
 Let $T$ a Hilbert-Schmidt operator on $\overline{H}^2(\Omega)\oplus L^2(\Omega)$. We assume that there exists $0\le \theta_1<\theta_2<\cdots<\theta_N<2\pi$ such that $\theta_k-\theta_{k-1}<\pi/2$ for $k=2,\cdots, N$ and $2\pi -\theta_N+\theta_1<\pi/2$ satisfying there exist $r_0>0$,  $C>0$ such that 
$\sup_{r\ge r_0}\| T_{re^{i\theta_k}} \|_0\le C$, for $k=1,\cdots, N$.  Moreover we assume there exists $(\lambda_j)$ such that $|\lambda_j|\to +\infty$ and
for all $f$ and $g$ in  $ \overline{H}^2(\Omega)\oplus L^2(\Omega)$, $(T_{\lambda_j}f|g)\to 0$.
Then the space spanned by the non zero generalized eigenfunctions of $T$ is dense in the adherence of the range of $T$.
\end{proposition}

\begin{prooff}
 As in Agmon~\cite[page 284]{Ag65} we define $F(\lambda)=(T_\lambda f|g)$ where $g$ is orthogonal to the generalized eigenfunctions. The goal is to prove that $F(0)=0$. As in Agmon we can prove that $F(\lambda)$ is analytic in $\C$ and bounded. Then $F(\lambda) $ is constant by Liouville theorem and as $(T_{\lambda_j}f|g)\to 0$ this implies that $F(\lambda)=0$.
\end{prooff}

\begin{prooff}[Proof of Theorem~\ref{Spect infini}] 

Before the proof we give some results on the links between the spectral decomposition of $S$ and $S^p$, where $S$ is an bounded operator on $\overline{H}^2(\Omega)\oplus L^2(\Omega) $.

 Let $\omega_j$ for $j=1,\cdots, p$, the roots of $z^p=1$. We have $z^p-1=\prod_{j=1}^p(z-\omega_j)$, in particular for $z=0$ we have $-1=\prod_{j=1}^p(-\omega_j)$. Thus we have
\begin{align}
 z^p-1=\prod_{j=1}^p(z-\omega_j)=\prod_{j=1}^p(-\omega_j)\prod_{j=1}^p(1-\omega_j^{-1}z)=-\prod_{j=1}^p(1-\omega_jz),\label{zp moins 1}
\end{align}
as we have $\{\omega_j,\  j=1,\cdots, p\}=\{\omega_j^{-1},\  j=1,\cdots, p\}$.

Applying~\eqref{zp moins 1} to $zS$,  we obtain
\begin{equation}\label{pd Rp}
( 1-z^pS^p)=\prod_{j=1}^p(1-\omega_jzS).
\end{equation}
If $(I-\omega_jzS)$ is invertible for all $j$, this implies that $( I-z^pS^p)$ is invertible. If for a fixed $j$, $(I-\omega_jzS)$ is not invertible, either $\Ker (I-\omega_jzS)\not= \{0\}$ this implies $\Ker (I-z^pS^p)\not= \{0\}$  or the range is not $\overline{H}^2(\Omega)\oplus L^2(\Omega)$ this implies that the range of $( I-z^pS^p)$ is not $\overline{H}^2(\Omega)\oplus L^2(\Omega)$. We deduce that 
$$
I-\omega_jzS \text{ is invertible for all } j \Leftrightarrow  I-z^pS^p \text{ is invertible}.
$$
If $ S^p$ is compact and $  I-z^pS^p $ is not invertible then by the Riesz theorem $z^{-p}$ is an eigenvalue of $S^p$ and there exists $k$ such that $\Ker (I-z^pS^p)^{k-1}\not=\Ker (I-z^pS^p)^{k}=\Ker (I-z^pS^p)^{k+1}$ and the dimension of $ \Ker (I-z^pS^p)^{k}$ is finite.

We will prove that all the eigenvalues of $S$ have the form $\omega_j z^{-1}$. Indeed $S$ is a operator on 
$\Ker (I-z^pS^p)^{k}$, then $S$ admits a spectral decomposition on $\Ker (I-z^pS^p)^{k}$. Let $u\not=0$ and $\lambda$ such that $u=\lambda Su$ then $z^pu=(\lambda zS)^pu$ thus $\lambda^{kp}(I-z^pS^p)^ku=(\lambda^p-z^p)^ku$ and $\lambda^p=z^p$ this implies $\lambda=\omega_j z$.

Now we prove that $\Ker (I-\omega_jzS)^k=\Ker (I-\omega_jzS)^{k+1}$. 
From \eqref{pd Rp}, we have
 $$u\in \Ker (I-\omega_jzS)^{k+1}\subset \Ker (I-z^pS^p)^{k+1}= \Ker (I-z^pS^p)^{k}.$$
We have
\begin{align*}
 (I-z^pS^p)^k&=\left( I-\left( I-(I-z\omega_jS)
\right)^p 
\right) ^k 
= \left( p(I-z\omega_jS)+\sum_{\mu=2}^p C_\mu(I-z\omega_jS)^\mu
\right) ^k\\
&=p^k(I-z\omega_jS)^k\left( I+\sum_{\mu=1}^{p-1}C'_\mu(I-z\omega_jS)^\mu
\right) ^k,
\end{align*}
 this implies $(I-z^pS^p)^ku=p^k(I-z\omega_jS)^ku=0$, which is the claim.

Obviously we can find the spectral decomposition of $S^p$ from the one of $S$. This proves that there exists $j$ such that $\Ker (I-\omega_jzS)^{k-1}\not=\Ker (I-\omega_jzS)^{k}=\Ker (I-\omega_jzS)^{k+1}$.

To prove the Theorem~\ref{Spect infini}, we fix $z$ as in the proof of Theorem~\ref{th : Existence}, we denote by $S=R_z$ and  we apply Proposition~\ref{prop infinie vp} to $T=S^p=R_{z}^p$, By Theorem~\ref{th est spect} and Lemma~\ref{lem est H-S}, $T: \overline H^2(\Omega)\oplus L^2 (\Omega)\to H^{2+2p}(\Omega)\oplus H^{2p}(\Omega)$ then $T$ is an Hilbert-Schmidt operator as $p>n/4$. We remark that 
\begin{align}\label{lien RS}
 S_\lambda=(R_{z})_\lambda=R_{z}(I-\lambda R_{z})^{-1}=(R_{z}^{-1}-\lambda)^{-1}=(B_{z}-\lambda)^{-1}=(B_0-z-\lambda)^{-1}=R_{z+\lambda}.
\end{align}
As $C_e$ is a closed cone, if $re^{i\theta}$ is not in $C_e$ for all $r$ large enough, then $z+re^{i\theta}$  is not in $C_e$ for all $r$ large enough (see \eqref{def C} for the notation $C_e$).

Using the previous remark and the estimate~\eqref{Est norm R} we have that $\|S_{re^{i\theta}} \|=\| R_{z+re^{i\theta}}\|_0$ is bounded uniformly with respect $r$ large enough, if $\theta\not=0$ and $re^{i\theta}\not\in C_e$. 

We prove the following formula
\begin{equation}\label{Lien R T}
 pz^{p-1}T_{z^p}=\sum_{k=1}^p\omega_k S_{\omega_kz}.
\end{equation}
Indeed we take the inverse of \eqref{pd Rp} when the formula make sense, we have
\begin{equation}\label{pd -1}
( 1-z^pS^p)^{-1}=\prod_{j=1}^p(1-\omega_jzS)^{-1}.
\end{equation}
 Derivate~\eqref{zp moins 1} with respect $z$, we obtain
\begin{equation*}
 pz^{p-1}=\sum_{k=1}^p\omega_k\prod_{j=1, j\not=k}^p(1-\omega_jz).
\end{equation*}
We apply this formula to $zS$, we obtain
\begin{equation}\label{R p-1}
 pz^{p-1}S^{p-1}=\sum_{k=1}^p\omega_k\prod_{j=1, j\not=k}^p(1-\omega_jzS).
\end{equation}

We multiply term by term \eqref{R p-1} and \eqref{pd -1}, multiplying the result by $S$, we obtain~\eqref{Lien R T}. 

If $S_{\omega_kz}$ is bounded uniformly for $r$ large enough and for all $k$, by \eqref{Lien R T} we have 
$\| T_{z^p} \|_0\le \frac{C}{|z|^{p-1}}$.  If we assume that $C_e$ is contained in a sector less than $\theta$ with $\theta<\pi/2$ and $\theta <2\pi/p$, the union of   $C_e$ and $C_e$ rotated  by angle $2k\pi/p$ does not give $\C$ and the Formula~\eqref{Lien R T} proves that we can find   the $\theta_j$'s satisfying the assumption of Proposition~\ref{prop infinie vp}. If $p\ge2$ the estimate on  $\| T_{z^p} \|_0 $ is stronger than the weak convergence. In case $p=1$, we have by the Theorem~\ref{th est spect} and Proposition~\ref{est v L2} with the notation $f=(f_1,f_2)$,
\begin{equation*}
 \| R_{z} f\|_{\overline{H}^2(\Omega)\oplus L^2(\Omega)}\le C\|\chi_\delta f_1 \|_{L^2(\Omega)}+C_\delta |z|^{-1}\|  f_1\|_{\overline{H}^2(\Omega)}+C  |z|^{-2}\| f_2 \|_{L^2(\Omega)},
\end{equation*}
if $|z|$ is large enough and $z\not\in C_e$. For $f_1$ and $\eps>0$ fixed we can choose $\delta>0$ such that $ C\|\chi_\delta f_1 \|_{L^2(\Omega)}\le \eps$. Then it is easy to prove that if $z_j$ is on a line such that $z_j=r_je^{i\theta}$ with $r_j\to + \infty$ , we have $\limsup \| R_{z_j} f\|_{\overline{H}^2(\Omega)\oplus L^2(\Omega)}\le \eps$. This prove that $ \| R_{z_j} f\|_{\overline{H}^2(\Omega)\oplus L^2(\Omega)}\to 0$ thus $ \| S_{z_j-z} f\|_{\overline{H}^2(\Omega)\oplus L^2(\Omega)}\to 0$. 

Now we prove that the adherence of $R_z(\overline H^2(\Omega)\oplus L^2(\Omega))$ is $H^2_0(\Omega)\oplus\{ v\in L^2(\Omega),\ \Delta v\in L^2(\Omega)\}$. Let $(u,v)\in H^2_0(\Omega)\oplus\{ v\in L^2(\Omega),\ \Delta v\in L^2(\Omega)\}$ we have $B_z(u,v)=(f,g)\in L^2(\Omega)\oplus  L^2(\Omega)$. Let $(f_n,g_n)\in  \overline H^2(\Omega)\oplus L^2(\Omega)$ such that $(f_n,g_n) \to (f,g)$ in $L^2(\Omega)\oplus  L^2(\Omega)$. We can take for instance $f_n$ and $g_n$ in $\Con_0^\infty(\Omega)$. We have by continuity $R_z(f_n,g_n)\to R_z(f,g)=(u,v)$ in $H^2_0(\Omega)\oplus\{ v\in L^2(\Omega),\ \Delta v\in L^2(\Omega)\}$ with the norm defined by  $\| v\|_{L^2(\Omega)}+\| \Delta v\|_{L^2(\Omega)}$ on $\{ v\in L^2(\Omega),\ \Delta v\in L^2(\Omega)\}$ and the usual norm on $H^2_0(\Omega)$.
\end{prooff}

\begin{prooff}[Proof of Theorem~\ref{Fct compt}]
Using \eqref{lien RS} as $z$ is fixed, to estimate the number of eigenvalues less than $t^2$ is equivalent to estimate the number of eigenvalues $\lambda$ such that $z+\lambda$  is less than $t^2$. That is in the sequel we estimate the number of $\lambda$ less than $t^2$ such that $\lambda^{-1}$ is a eigenvalue of $S=R_z$. 

 We have shown in the proof of Theorem~\ref{Spect infini} that $\| T_{z^p} \|_0\le \frac{C}{|z|^{p-1}}$ if  $\omega_kz$ is on a line $d(r_0,\theta)=\{z\in \C,\ z=re^{i\theta}, r\ge r_0\}\subset C_e$. As $(1-z^p T)^{-1}=1+z^p T_{z^p}$, we obtain for $|z|\ge 1$, $\| (1-z^p T)^{-1} \|_0\le 1+|z|^p\| T_{z^p} \|_0\le C|z|$.

We have, by  $ \| T_{z^p}\|_{2p}\le \|  T\|_{2p}\| (1-z^p T)^{-1}\|_0\le C|z|$. We obtain  from Lemma~\ref{lem est H-S}
$$ \tno{T_{z^p}}\le C\| T_{z^p}\| _{2p}^{n/(4p)}   \| T_{z^p}\|_0^{1-n/(4p)}\le C|z|^{1-p+n/4}.$$

We remind \cite[Theorem 12.14]{Ag65} if $T$ is Hilbert-Schmidt, we have $\sum |\mu_i|^2\le \tno{T}^2$ where $\mu_j \not=0$ are the eigenvalues counted with multiplicities.

Let $\lambda_j$ such that $\lambda_j^{-1}$ is a eigenvalue of $S$, we find that $\frac{1}{\lambda_j^p-z^p} $ is a eigenvalue of $T_{z^p}$. 
We obtain 
$$
\sum_j\frac{1}{|\lambda_j^p-z^p|^2} \le\tno{T_{z^p}}^2\le C |z|^{2-2p+n/2}.
$$
If $|\lambda_j|\le t^2$ and taking $z \in d(r_0,\theta)$ satisfying $|z|=t^2$, we have $ |\lambda_j^p-z^p|\le 2t^{2p}$. Then we have
$$
\sum_{|\lambda_j|\le t^2}\frac 1{4t^{4p}}\le \sum_j\frac{1}{|\lambda^p_j-z^p|^2} \le\tno{T_{z^p}}^2\le Ct^{4-4p+n}.
$$
Then we obtain that $N(t)\le Ct^{4+n}$.
\end{prooff}

\section{Estimate on the resolvent}\label{sec est res}

\subsection{Upper bound}\label{sub upper}

In this section we prove the Theorem~\ref{Maj Res}.
We recall the well-known Green's formula. For regular functions $u$ and $v$, we have
$$
\int_\Omega(u\Delta v-v\Delta u)dx=\int_{\d \Omega}(u\d_\nu v-v\d_\nu u)ds,
$$
where $\d_\nu$ is the exterior normal derivative on $\d \Omega$ and $ds$ is the surface measure on $\d \Omega$. Here we work with smooth functions. As the problem is well-posed by Theorem~\ref{th : Existence} we can apply the estimate for non smooth functions by passing to the limit in the estimate.

\begin{equation}\label{eq : n complex}
 \left\{
	\begin{array}{l}
	 	\Delta w +k^2n(x) w=f \text{ in } \Omega,\\ 
		\Delta v +k^2 v=g \text{ in } \Omega,\\
		w=v \text{ on } \d \Omega,\\
		\d_\nu w=\d_\nu v  \text{ on } \d \Omega.
	\end{array}\right.
\end{equation}

By Green's formula and  \eqref{eq : n complex} we have

\begin{align*}
 \int_\Omega(w\Delta \bar w-\bar w\Delta w)dx &=\int_{\d_\Omega}(w\d_\nu \bar w-\bar w\d_\nu w)ds\\
&=\int_\Omega [w(\bar f-k^2\bar n\bar w)-\bar w(f-k^2nw)]dx\\
&=\int_\Omega [w\bar f-\bar w f-2ik^2\Im n |w|^2]dx,
\end{align*}
and

\begin{align*}
 \int_\Omega(v\Delta \bar v-\bar v\Delta v)dx &=\int_{\d_\Omega}(v\d_\nu \bar v-\bar v\d_\nu v)ds\\
&=\int_\Omega [v(\bar g-k^2\bar v)-\bar v(g-k^2v)]dx\\
&=\int_\Omega[ v\bar g-\bar v g ]dx.
\end{align*}

Using the boundary condition in \eqref{eq : n complex}, we obtain,

$$
\int_\Omega (w\bar f-\bar wf)dx-\int_\Omega (v\bar g-\bar vg)dx=2i\int _\Omega k^2|w|^2\Im ndx.
$$
Thus we deduce,
\begin{equation}\label{est : Im n positif}
 \delta\int_\omega k^2|w|^2\le \|v\|\|g\|+\| w\|\| f\|,
\end{equation}
where $\omega=\{x\in \Omega,\ \Im n(x)\ge\delta\}$.
\begin{remark}
 In the case where $n=n_1+in_2/k$ we have $\Im n=n_2/k$, and in the previous computations we must change the left hand side of~\eqref{est : Im n positif} by $ \delta\int_\omega k|w|^2 $ where $\omega=\{x\in \Omega,\  n_2(x)\ge\delta\}$. We let to the reader to check that the rest of the proof does not change with this new estimate. Indeed the powers of $k$ do not play any role with respect the estimates by $e^{Ck}$.
\end{remark}

We recall the interpolation estimate. We can find this type of estimate in \cite[ Section  3, Formulas (1) and (2)]{LR95}, \cite[Theorem 3]{LR97}, \cite[Proposition 1.2]{CR12}. The estimate~\eqref{int est 1} does not appear in this literature, but we can prove it following the same ways. Indeed in the Carleman estimate used to prove the interpolation estimates, we estimates also the boundary terms but in the previous mentioned paper we did not need the boundary term in the interpolation estimates.

 Let $X=(-3,3)\times \Omega$, $Y=(-2,2)\times \Omega$, and ${\cal O}=(-1,1)\times \omega$. We denote by $\d Y=(-2,2)\times\d \Omega$. Then there exist $\delta>0$ and $C>0$ such that for all $W\in \overline H^1(X)$ such that $\d_s^2W+\Delta W\in L^2(X)$, $W_{|\partial Y}\in H^1(\partial Y)$, $\d_\nu W_{|\partial Y}\in L^2(\d Y)$ , we have
\begin{align}
&\|W\|_{\overline H^1(Y)}\!\!+| W_{|\partial Y}|_{H^1(\partial Y)}\!+| \d_\nu W_{|\partial Y}|_{L^2(\d Y)}\le\! C
\left( 
\| \d_s^2W\!+\Delta W\|_{L^2(X)}\!+\| W\|_{L^2({\cal O})}
\right) ^\delta \!
\|W\|_{\overline H^1(X)}^{1-\delta} \label{int est 1}
\\
&\|W\|_{\overline H^1(Y)}\le\! C
\left( 
\| \d_s^2W+\Delta W\|_{L^2(X)}+| W_{|\partial Y}|_{H^1(\partial Y)}+| \d_\nu W_{|\partial Y}|_{L^2(\d Y)}\right) 
^\delta \!
\| W\|_{\overline H^1(X)}^{1-\delta},\label{int est 2}
\end{align}
where $s$ is an additional variable. This variable allows us to give an estimate uniform with respect the large parameter $k$. We shall see that in the sequel.

Let $W(s,x)= e^{sk}w(x)$ where $\Delta w+k^2w=f$ in $\Omega$. We have $\d_s^2 W+\Delta W=e^{sk}f$ and we can obtain the following estimates for a $C>0$,
\begin{align*}
 &\| w\|_{\overline H^1(\Omega)}\le C \|  W\|_{\overline H^1(Y)},\\
&| w_{|\d \Omega}|_{H^1(\d \Omega)}\le C | W_{|\d Y}|_{\overline H^1(\d Y)},\\
&|\d_\nu w_{|\d \Omega}|_{L^2(\d \Omega)}\le C| \d _\nu W_{|\d Y}|_{L^2(\d Y)},\\
&\| \d_s^2W+\Delta W\|_{L^2(X)}\le Ce^{3k}\| f\|_{L^2(\Omega)},\\
&\| W\|_{L^2({\cal O})}\le C e^{k}\| w\|_{L^2(\omega)},\\
&\| W\| _{\overline H^1(X)}\le Ce^{4k}\| w\| _{\overline H^1(\Omega)}.
\end{align*}
By the interpolation estimate~\eqref{int est 1},
there exists $C>0$ such that for all $w\in\overline H^1(\Omega)$, satisfying $\d_\nu  w_{|\d \Omega}\in L^2(\d \Omega) $ and $ w_{| \d \Omega}\in  H^1(\d \Omega)$ solution of $\Delta w+k^2w=f$ in $\Omega$, we have
\begin{equation}\label{est. exp.}
| (\d_\nu  w)_{|\d \Omega} |_{L^2(\d \Omega)} +| w_{| \d \Omega} |_{H^1(\d \Omega)}+ \| w \|_{\overline H^1(\Omega)}\le Ce^{Ck}(\| w \|_{L^2(\omega)}+\| f \|).
\end{equation}
Using \eqref{est : Im n positif}, \eqref{est. exp.} and $Ce^{Ck}\| w\|\| f\|\le (1/2)\| w \|^2+ C^2e^{2Ck}\|  f\| ^2$,  we have for a $C>0$
\begin{equation}\label{est. w}
 \| w \|^2_{\overline H^1(\Omega)}\le Ce^{Ck}(\| f \| ^2+\| v\| \| g\|).
\end{equation}
Following the same way, we denote by $W(s,x)=e^{sk}v(x)$ and we apply the interpolation estimate~\eqref{int est 2}, we obtain on 
 $v$  the estimate
\begin{equation*}
 \| v \|_{\overline H^1(\Omega)}\le Ce^{Ck}(\| g \|+| v_{|\d \Omega}|_{H^1(\d \Omega)}+|  \d_\nu v_{|\d \Omega}|_{L^2(\d \Omega)}).
\end{equation*}
Taking account the boundary condition in \eqref{eq : n complex} and \eqref{est. exp.} we have
\begin{align}
  \| v \|_{\overline H^1(\Omega)}^2&\le Ce^{Ck}(\| g \|^2+| w_{|\d \Omega}|_{H^1(\d \Omega)}^2+|  \d_\nu w_{\d \Omega}|_{L^2(\d \Omega)}^2)  \notag \\
& \le Ce^{Ck}(\| g \|^2+\| f \|^2 + \| w \|_{L^2(\omega)}^2)  \notag\\
&\le Ce^{Ck}(\| g \|^2+\| f \|^2 +\|v\|\|g\|+\| w\|\| f\|) \text{ by \eqref{est : Im n positif}}.\notag
\end{align}
This estimate and \eqref{est. w} give
\begin{equation*}
 \| v \|_{\overline H^1(\Omega)}+\| w\|_{\overline H^1(\Omega)}\le Ce^{Ck}(\| f\| +\| g \|).
\end{equation*}
This implies the estimate on the $L^2$ norm on $v$ and $w$, which gives the Theorem~\ref{Maj Res} by density.

\subsection{Lower bound}\label{sub lower}

Here we use the results proved first by Davies~\cite{Da99} in one dimension,  by Zworski~\cite{Zw01} for Schrödinger operators in $n$ dimension and by  Dencker, Sjöstrand and Zworski~\cite{DSZ04} for more general sub-elliptic operators. This allows  to obtain a lower bound on the resolvent. 

We recall here the theorem given by  Dencker, Sjöstrand and Zworski.
\begin{theorem} [Theorem 1.1 \cite{DSZ04}]  \label{Th DSZ}
Let  $V\in\Con^\infty(\R^n)$. Then, for any $z\in \{ \xi^2+V(x),\ (x,\xi)\in\R^n$, $\Im \est{\xi,\d_xV(x)}\not=0\}$, there exists $h_0>0$ such that for all $h\in(0,h_0)$, there exists $u(h)\in L^2(\R^n)$ with the property
$$
\| (-h^2\Delta+V(x)-z)u(h)\|={\cal O}(h^\infty)\| u(h) \|.
$$
In addition, $u(h)$ is localized to a point $(x_0,\xi_0)$ in phase  with $\xi_0^2+V(x_0)=z$.

More precisely, $WF_h(u)=\{ (x_0,\xi_0 )\}$, where $WF_h(u)$ is  the semi-classical wave front set.

If the potential is real analytic, then we can replace $h^\infty$ by $\exp(-1/Ch)$.
\end{theorem}

A consequence of the microlocal  localization of $u$, we can cut-off $u$ such that its support is in a \nhd of $x_0$. The Theorem~\ref{Th DSZ} implies that if $z$ is in the resolvent set, $\|(-h^2\Delta+V(x)-z)^{-1} \|\ge C_Nh^{-N}$ for all $N$, in $\Con^\infty $ case and $\|(-h^2\Delta+V(x)-z)^{-1} \|\ge C\exp(C/h)$ in analytic case.

\begin{prooff}[Proof of Theorem~\ref{min Res}]
 We set $z_0=-(1+m(x_0))^{-1}|\xi_0|^2$, we  set $V(x)=z_0(1+m(x))$. We have $|\xi_0|^2+V(x_0)=0$ and 
 \begin{align*}
 \Im (\xi_0\d_x V(x_0))&=\Im ( z_0 \xi_0\d_xm(x_0))=-|1+m(x_0)|^{-2}|\xi_0|^2\Im ((1+\overline{ m}(x_0))\xi_0\d_xm(x_0))\\
  &= -|1+m(x_0)|^{-2}|\xi_0|^2\Im (\bar n(x_0)\xi_0\d_xn(x_0)) \not=0,
 \end{align*}

 by assumption. By Theorem~\ref{Th DSZ} there exists $u(h)$ such that $\|(-h^2\Delta +V(x))u(h)\|={\cal O}(h^\infty)\| u(h) \|$ or $={\cal O}(e^{-C/h)}\| u(h) \|$ if $m$ is analytic. Define  by $f=\Delta u(h)-h^{-2}z_0(1+m(x))u(h)$, we have $\tilde R_{k^2}(f,0)=(u(h),0)$ with $k^2=-h^{-2}z_0$. We remark that $u(h)$ is localized in a \nhd of $x_0$, in particular $u(h)$ is null in a \nhd of $\d\Omega$, then $(u(h),0)$ satisfies the boundary conditions. This implies the Theorem~\ref{min Res}.
\end{prooff}

\appendix

\section{Notation and recall on pseudo-differential calculus}\label{App pseudo}

\subsection{Sobolev spaces and pseudo-differential operators}

We introduce some notation for the Sobolev spaces.

We denote the semi-classical $H^s$ norm by $\| w\|_{H^s_{sc}}^2=\int (1+h^2|\xi|^2)^s|\hat u(\xi)|^2d\xi$. On a compact manifold we define the semi-classical $H^s$ using local coordinates. To distinguish norm on spaces of dimension $n$ and  dimension $n-1$, we denote the semi-classical $H^s$ norm on $\partial \Omega$ by $| \cdot |_{H^s_{sc}(\d \Omega)}$.  Let $w$ a distribution on $\Omega$, we denote by $\| w\|_{\overline{ H}^s_{sc}(\Omega)}=\inf\{ \| \beta \|_{H^s_{sc}}, \text{ where }\beta_{|\Omega}=w  \}$. We recall that we denote by $D=(h/i)\d$,  and if $s$ is an integer the quantity $\sum_{|\alpha|\le s}\| D^\alpha w \|^2_{L^2(\Omega)}$ is equivalent uniformly with respect $h$ to $\| w\|^2_{\overline{ H}^s_{sc}(\Omega)}$. 

In the context of semi-classical $H^s$ space we have the following trace formula, for $s> 0$, for all $w\in \overline H^{s+1/2}_{sc}(\Omega)$ we have
$$
| w_{|\d \Omega}|_{H^s_{sc}(\d \Omega)}\lesssim h^{-1/2}\| w\|_{\overline H^{s+1/2}_{sc}(\Omega)},
$$
where $w_{|\d\Omega}(x_0)$ means the limit of $w(x)$ when $x\in \Omega$  goes to $x_0$.

We recall the pseudo-differential tools.
Let $a(x,\xi)$ in  $\Con^\infty(\R^n\times\R^n)$ we say that $a$ is a symbol of order $m$ if for all $\alpha, \beta\in\N^n$, there exist $C_{\alpha, \beta}>0$,  such that
$$
|\d_x^\alpha\d_\xi^\beta a(x,\xi)|\le C_{\alpha, \beta}\est{\xi}^{m-|\beta|},
$$
where $\est{\xi}^2=1+|\xi|^2$. In particular a polynomial in $\xi$ of order $m$ with coefficients in $\Con^\infty(\R^n)$ with all bounded derivatives, is a symbol of order $m$.

To a symbol we can associate an semi-classical operator by the following formula
$$
\Op(a)u=a(x,D)u= \frac{1}{(2\pi)^n}\int e^{ix\xi}a(x,h\xi)\hat u(\xi)d\xi= \frac{1}{(2h\pi)^n}\int e^{ix\xi/h}a(x,\xi)\hat u(\xi/h)d\xi.
$$
This formula makes sense for $u\in\S(\R^n)$ and we can extend to $u\in H^s$. For $a$, a symbol of order $m$, there exists $C>0$ such that for all $u\in H^s$,
$$
\| a(x,D)u\|_{H_{sc}^{s-m}}\le C\| u\|_{H_{sc}^{s}}.
$$

We can compose the pseudo-differential operators, and we have the following result. Let $a$ a symbol of order $m$ and $b$ a symbol of order $k$, there exists $c$ a symbol of order $m+k$ such that $a(x,D)\circ b(x,D)=c(x,D)$. Moreover there exists a symbol $d$ of order $m+k-1$ such that $c(x,D)= (ab)(x,D)+hd(x,D)$. This means that up $h$ the composition of two operators is the operators associated with the product of symbols.

We can inverse the elliptic symbol, more precisely, let $a$ a symbol of order $m$ satisfying there exists $C>0$ such that for all $(x,\xi)\in\R^n\times\R^n$ we have $|a(x,\xi)|\ge C\est{\xi}^m$. Then for all $N>0$ there exist $b$ a symbol of order $-m$ and $r$ a symbol of order $-N$ such that $b(x,D)\circ a(x,D)=I+hr(x,D)$.
We can localized this result. Let $K $ a closed set of $\R^n$, we assume that there exists $C>0$ such that for all $(x,\xi)\in K\times\R^n$, $|a(x,\xi)|\ge C\est{\xi}^m$ for all $\chi\in\Con^\infty_0(\R^n)$ supported in $K$, there exists $b$ a symbol of order $-m$ and $r$ a symbol of order $-N$ such that $b(x,D)\circ a(x,D)=\chi(x)+hr(x,D)$. In  both cases we say that $b$ is a parametrix for $a$.

We can also define pseudo-differential on a smooth compact  manifold without boundary. We shall use freely the result on $\R^n$ in the context of manifolds. To distinguish both cases we denote by $\Op a$ the operators on $\R^n$ and by $\op a$ the operators on a manifold of dimension $n-1$ or on $\R^{n-1}\cong \{ x\in\R^n, \ x_n=0\}$.

We use also spaces $\overline H^s_{sc}$ and the pseudo-differential calculus on these spaces. In general that requires introduction of the delicate notion of ``transmission condition'' (see Boutet de Monvel~\cite{Bo71}), to avoid that we follow the Hörmander's strategy (see~\cite[Appendix B]{Hobo83}) adapted for the parametrices which are particular cases of operators satisfying the ``transmission condition''. We recall some results proved by Hörmander in the context of classical $\overline H^s$ spaces. The adaptation to the $\overline H^s_{sc}$ spaces is easy and we give here only the results and some ideas of proof. Here we give the result in a half space $\R^n_+=\{ x\in\R^n, \ x_n>0\}$.
For simplicity we denote by $\overline H^s_{sc}=\overline H^s_{sc}(\R^n_+)$.
In the proof we need introduce a space $\overline H^{m,s}_{sc}$. First we say that $u\in H_{sc}^{m,s}(\R^n)=H_{sc}^{m,s}$ if $\| u\|^2_{H_{sc}^{m,s}}=\int \est{h\xi}^{2m}\est{h\xi'}^{2s}|u(\xi)|^2d\xi<\infty$ where $\xi=(\xi',\xi_n)$. As for the $H^s$ space we say that $ u\in \overline H^{m,s}_{sc}$ where $u$ is a distribution in $\D'(\R^n_+)$ 
 if there exists $v\in  H^{m,s}_{sc}$ such that $u=v_{|x_n>0}$ and we denote $\| u\|_{\overline H^{m,s}_{sc}}=\inf\{ \| v\|_{H^{m,s}_{sc}}, \text{ where }{v\in H^{m,s}, \text{ such that }v_{|x_n>0}=u}\}$.

We can easily see that if $u\in \overline H_{sc}^{0,s}$ then $\underline u\in H_{sc}^{0,s}$ and $\| u \|_{ \overline H_{sc}^{0,s}}=\|\underline u \|_{  H_{sc}^{0,s}}$, (see \eqref{prolongement 0} for the definition of $\underline u$).

We can extend the Theorem~B.2.3 given by Hörmander~\cite{Hobo83}
\begin{align}
u\in \overline H^{m,s}_{sc}&\Leftrightarrow D_{x_n}u\in \overline H^{m-1,s}_{sc}\text{ and } u\in \overline H^{m-1,s+1}_{sc}\notag \\
&\Leftrightarrow  
 D_{x_n}^2u\in \overline H^{m-2,s}_{sc}\text{, }D_{x_n}u\in \overline H^{m-2,s+1}_{sc}
  \text{ and }u\in \overline H^{m-2,s+2}_{sc}.\label{Hm,s}
\end{align}
Of course the natural norms on these spaces are equivalent.

We can use the Theorem~B.2.9 from \cite{Hobo83} in the following form adapted to our context. We denote by $P=D_{x_n}^2+R(x,D')+\alpha(x)D_{x_n}$ a differential operator of second order. We have for all $k\in \R$,

\begin{equation}\label{car Hm,s}
 u\in  \overline H^{m-k,s+k}_{sc}\text { and }  Pu\in \overline H^{m-2,s}_{sc} \Rightarrow u\in \overline H^{m,s}_{sc}.
\end{equation}
Of course this is trivial if $k\le 0$. For $k>0$ we can prove this by recurrence on $k$. The idea is the following, we can write $D_{x_n}^2u=Pu-R(x,D')u-\alpha D_{x_n}u$ and for $k=1$ this formula
 implies that $D^2_{x_n}u\in \overline H_{sc}^{m-2,s}$. Then we can apply \eqref{Hm,s} to obtain  the result in this case. The recurrence is  easy.

The previous results are useful to prove that the parametrix of a elliptic operator is a mapping  on the $\overline H^s $ space. 

Let $P$ a differential operator of second order, elliptic, i.e.  there exists $C>0$ such that $\forall x\in\R^n$, $\forall \xi\in\R^n$, $|p(x,\xi)|\ge C\est{\xi}^2$, and let $Q$ a parametrix such that $QP=Id +hK$ where $K$ is of order $-N$ where $N>0$. For a distribution $w$  in $\R^n$, we denote  by $rw=w_{|x_n>0}$. The action of $Q$ on $\overline H^s_{sc}$ is given by the  formula $rQ\underline u$ which make sense if $\underline u$ make sense (it is the case for if $u\in L^2(\Omega)$). We have the following result.
If $s\in[0,N]$, there exists $C>0$ such that for all $u\in\overline H^s_{sc}$, we have 

\begin{equation}\label{eq ope para}
\| rQ\underline u\|_{\overline H^{s+2}_{sc}}\le C\| u \|_{\overline H^{s}_{sc}} . 
\end{equation}
Remark, here $\underline u$ is at least in $L^2$ so $rQ\underline u$ make sense.
First we prove that $r Q\underline u\in \overline H^{s+2-k,k}_{sc}$ where $k\ge s$. It is enough to prove that for all $\alpha\in\N^{n-1}$, $|\alpha|\le k$, $D'^\alpha Qu\in \overline H^{s+2-k}_{sc}$. By pseudo-differential calculus we have $D'^\alpha Q=\sum_{|\beta|\le |\alpha|}Q_\beta D'^\beta$, where $Q_\beta$ is of order $-2$. We have 

\begin{align*}
\| rD'^\alpha Q\underline u\|_ {\overline H^{s+2-k}_{sc}}&\le \sum_{|\beta|\le |\alpha|}\| Q_\beta D'^\beta \underline u \|_ { H^{s+2-k}_{sc}}\le C  \sum_{|\beta|\le |\alpha|}\|  D'^\beta \underline u \|_ { H^{s-k}_{sc}}\le C\| \underline u\|_{ H^{s-k,k}_{sc}}\\
&\le   C\| \underline u\|_{ H^{0,s}_{sc}}\le  C\| u \|_{\overline H^{0,s}_{sc}}\le C\| u \|_{\overline H^{s}_{sc}}.
\end{align*}

It is well-known that we have also $PQ=Id+h\tilde K$ where $\tilde K$ is of order $-N$. Indeed there exists a $\tilde Q$ such that $P\tilde Q=Id+h\tilde K$ where $\tilde K$ is of order $-N$ and it is easy to prove that $Q=\tilde Q+hK'$ where $K'$ is of order $-N-2$. Thus we have, as $P$ is a differential operator, $PrQ\underline u=rPQ\underline u=r\underline u+hr\tilde K\underline u= u+hr\tilde K\underline u$. As $\| \tilde K\underline u\|_{H^N}\le C \| \underline u\|_{L^2}\le \|  u\|_{\overline H^s_{sc}}$, we obtain for $s\in[0,N]$, $ \| PrQ\underline u\|_{\overline H^s_{sc}}  \le C \|  u\|_{\overline H^s_{sc}}$. 
Then we have $rQ\underline u\in \overline H^{s+2-k,k}_{sc}$ and $PrQ\underline u\in \overline H^s_{sc}$ then from \eqref{car Hm,s} this implies $ \| rQ\underline u\|_{\overline H^{s+2}_{sc}}  \le C \|  u\|_{\overline H^s_{sc}}$.

We need also  regularity results for  $rQ(\gamma \otimes D^k_{x_n}\delta_{x_n=0})$, where $\gamma\in H^s_{sc}(\R^{n-1})$.
First we remark $\int h^{2k }|\xi_n|^{2k} (1+h^2|\xi_n|^2+|\xi'|^2)^\nu d\xi_n\le C\est{\xi'}^{2k+2\nu+1}/h$ if $\nu+k<-1/2$, then by direct computation we have for all $\gamma\in H^s(\R^{n-1})$, 
$\| \gamma \otimes D_{x_n}^k\delta_{x_n=0} \|_{H^{\nu-k-1/2,s-\nu}}\le \frac C{\sqrt h}|\gamma  |_{H^s(\R^{n-1})}$, if $\nu<k$. For $j\in\N$, we have if $s-k< j$,  following the same computation as for computing $rQ\underline u$,

$$
\|  Q( \gamma \otimes D_{x_n}^k\delta_{x_n=0})\|_{H^{s-j-k+3/2,j}}\le C\|   \gamma \otimes D_{x_n}^k\delta_{x_n=0}\|_{H^{s-j-k-1/2,j}}\le  \frac C{\sqrt h}|\gamma  |_{H^s(\R^{n-1})}.
$$

As $r(\gamma \otimes D_{x_n}^k\delta_{x_n=0})=0$ we have $PrQ(\gamma \otimes D_{x_n}^k\delta_{x_n=0})=   hrK(\gamma \otimes D_{x_n}^k\delta_{x_n=0})$. We deduce,
$$
\|PrQ(\gamma \otimes D_{x_n}^k\delta_{x_n=0}) \|_{\overline H^{s-k-1/2}_{sc}}\le C\| \gamma \otimes D_{x_n}^k\delta_{x_n=0}\|_{H^{s-k-N-1/2}_{sc}}\le \frac C{\sqrt h}|\gamma  |_{H^s(\R^{n-1})},
$$
if $s-k<N$. Thus from \eqref{car Hm,s} we obtain if $s-k<N$,

\begin{equation}\label{app est int par trace}
\| rQ(\gamma \otimes D_{x_n}^k\delta_{x_n=0})\|_{ \overline H^{s-k+3/2}(\R^n_+)}\le \frac C{\sqrt h}|\gamma  |_{H^s(\R^{n-1})}.
\end{equation}

\begin{remark}
 In the section~\ref{sub v bord} we apply the previous result to a local parametrix. Indeed we can construct the local parametrix with a global parametrix. We can extend $P$ to have a global elliptic operator $\tilde P$ such that $\tilde P=P $ in a domain $W$ where $P$ is elliptic.
Let $Q$ a parametrix of $\tilde P$ such that $Q\tilde P=Id+hK$ where $K$ is an operator of order $-N$.
Let $\chi_1$ and $\chi_2$ functions in $\Con^\infty$ compactly supported in $W$ such that $\chi_2=1$ on the support of $\chi_1$. By pseudo-differential calculus, we have $\chi_1 Q\chi_2=\chi_1Q+hK'$ where $K'$ is an operator of order $-N-2$. Then we have $\chi_1Q\tilde P=\chi_1+h\chi_1 K$ and $\chi_1Q\tilde P=\chi_1 Q\chi_2 \tilde P-hK'\tilde P$. As $\tilde P=P$ on the support of $\chi_2$ we have $\chi_1 Q\chi_2 P=\chi_1+hK''$ where $K''$ is an operator of order $-N$. Then $\chi_1Q\chi_2$ is a local parametrix of $P$. It is easy to see that we can replace $Q$ by $\chi_1Q\chi_2$ in~\eqref{eq ope para} and in~\eqref{app est int par trace}.
\end{remark}

\subsection{Properties on the roots and parametrices}\label{sb : roots param.}

We use some properties of the roots of $\xi_n^2+R(x,\xi')-\mu$ and $a(\xi_n^2+R(x,\xi'))-\mu$. By assumption these polynomials have not real roots and it is easy to see that for $\xi'$ large enough the imaginary parts have different signs. In particular the roots are simple thus smooth and the roots are symbols of order $1$. Actually, for instance for $\xi_n^2+R(x,\xi')-\mu$, (the proof for $a(\xi_n^2+R(x,\xi'))-\mu$ is similar and left to the reader) the roots have, for $|\xi'|$ large enough, the following form $\pm i\sqrt{R(x,\xi')}+z_\pm(\mu,1/\sqrt{R(x,\xi')})$, where $z_\pm$ is a solution to $z_\pm \mp i z_\pm^2s/{2}+i\mu s/2=0$ in a \nhd of $s=0$. This expression implies that the roots are symbols of order 1.

The parametrices used, denoted by $Q$ and $\tilde Q$ have a particular structure we give here.
The symbol of $P$ is a polynomial having the following form, $p_2+hp_1$ where $p_j$ are polynomial of degree $j$. We seek a parametrix with symbol given formally by $q=q_{-2}+hq_{-3}+h^2q_{-4}+\cdots$, where $q_j $ are symbol of order $j$. If we denote by $q\circ p$ the asymptotic expansion of the symbol of $\Op (q)\Op (p)$, we have,
$$
q\circ p=q_{-2}p_2+\sum \frac{h^{k-j+|\alpha|}}{\alpha ! i^{|\alpha|}}\d_\xi^\alpha q_{-k}\d_x^\alpha p_j,
$$
where in the sum we have $j=2$ or $1$, $k\ge 2$, $\alpha\in\N^n$ and, $|\alpha|\ge 1$ or $j=1$. In particular we have $k-j+|\alpha|\ge 1$. We choose $q_{-2}=1/p_2$, and to cancel the terms with the same power in $h$ we have
\begin{equation}\label{dev param}
q_{-\nu-2}=\frac{1}{p_2}\sum  \frac{1}{\alpha ! i^{|\alpha|}}\d_\xi^\alpha q_{-k}\d_x^\alpha p_j,
\end{equation}
where $\nu\ge 1$, $k-j+|\alpha|=\nu$, $j=2,1$, $k\ge 2$, $\alpha\in\N^n$, $|\alpha|\ge 1$ or $j=1$. In particular the sum is finite and $k\le \nu+1$.
We claim now that 
$$
q_{-\nu}=\frac{S_{3\nu-6}}{p_2^{2\nu-3}} \text{ for }\nu\ge2,
$$
where $S_\mu$ is a polynomial of degree $\mu$.

Clearly this is true for $\nu=2$. We verify that for $k\le \nu+1$,  $\d_\xi^\alpha q_{-k}=\frac{\tilde S_{3k-6+|\alpha|}}{p_2^{2k-3+|\alpha|}}$, where $\tilde S_\mu$ is a polynomial of degree $\mu$. The parameters satisfy $j\le 2$, $k-j+|\alpha|=\nu$ and $k\le \nu+1$, then the power of $p_2$ in \eqref{dev param} is $2k-2+|\alpha|=\nu+k+j-2\le 2 \nu+j-1\le 2\nu+1=2(\nu+2)-3$. The degree of the numerator is $3k-6+|\alpha|+2j=\nu+2k-6+2j\le 3\nu -4+2j\le 3\nu=3(\nu+2)-6$. This gives the claim.

We need to compute for $\gamma(x')$,

\begin{equation}\label{rest term bord}
 \left[
Q\left( 
\frac{h}{i}\gamma \otimes D^k_{x_n}\delta_{x_n=0}
\right) 
 \right]_{|x_n=0}=\frac{1}{(2h\pi)^{n-1}}\int e^{ix'\xi'/h}\tilde q(x',\xi')\hat \gamma(\xi'/h)d\xi',
 \end{equation}
where formally $\tilde q(x',\xi')=\left( 
\frac{1}{2i\pi}\int_\R e^{ix_n\xi_n/h}q(x,\xi)\xi_n^kd\xi_n
\right) _{|x_n=0}$. It is not clear that $\tilde q$ is well defined in general but in the following lemma we prove this is true if $q$ is a rational function, and in this case~\eqref{rest term bord} make sense.

\begin{lemma}\label{lem trace}
 Let $\nu\in\N^*$, let $S_{\nu}(x,\xi)$ a polynomial of order $\nu$ with respect $\xi_n $  and we assume that the coefficient of $\xi_n^j$ is a symbol in $\xi'$ of order $\nu-j$. Let $p$ a polynomial of degree $d$ in $\xi$ and a symbol of order $d$.
We assume that $p(x,\xi)=\xi_n^d+\sum_{j=0}^{d-1}\xi_n^ja_{d-j}(x,\xi')$ where $a_{d-j}$ are polynomials of order $d-j$. 
Moreover we assume there exists $\delta>0$ such that $\forall x\in\R^n$, $\forall \xi\in\R^n$, $|p(x,\xi)|\ge \delta \est{\xi}^d$.
Then
$$
\left(  \int_\R e^{i\frac{x_n}{h}\xi_n}\frac{S_{\nu}(x,\xi)}{p(x,\xi)}d\xi_n\right) _{|x_n=0},
$$
is a symbol of order $ \nu-d+1$. 
\end{lemma}

\begin{prooff}
The integral $\int_\R e^{i\frac{x_n}{h}\xi_n}\frac{S_{\nu}(x,\xi)}{p(x,\xi)}d\xi_n$ converges for $x_n>0$. For $(x,\xi')$ fixed, we can change the integration contour by $\Gamma =[-D\est{\xi'}D\est{\xi'}]\cup \{ z\in\C,\ |z|=D\est{\xi'},\ \Im z>0  \}$, where $D $ will be chosen later. Indeed the integral does not depend of $D$ if $D$ large enough and the integral on $\{ z\in\C,\ |z|=D\est{\xi'},\ \Im z>0  \}$ goes to $0$ if $D$ goes to $+\infty$. Now we integrate on a compact set and we can take the limit when $x_n$ goes to $0^+$. We have the  following quantity to control.
$$
A=\left(  \int_\Gamma\frac{S_{\nu}(0,x',\xi)}{p(0,x',\xi)}d\xi_n\right) .
$$
On $\Gamma$ we have $|S_{\nu}(0,x',\xi)|\le C\est{\xi'}^\nu$. For $ \xi_n\in\{ z\in\C,\ |z|=D\est{\xi'},\ \Im z>0  \}$ we have $|p(0,x',\xi)|\ge |\xi_n|^d-\sum_{j=0}^{d-1}|\xi_n|^j |a_{d-j}(0,x',\xi')|\ge D^d\est{\xi'}^d(1-dD^{-1}C)\ge \est{\xi'}^dD^d/2$, where we have $|a_{d-j}(0,x',\xi')|\le C\est{\xi'}^{d-j}$ and chosen $D$ such that $D\ge \max (1,2dC)$. 
Then by assumption for all $\xi_n\in\Gamma$ we have $|p(0,x',\xi)|\ge \delta'\est{\xi'}$.
As the length of $\gamma$ is less than $K\est{\xi'}$ we obtain $A\le K'\est{\xi'}^{\nu-d+1}$. We can obtain the estimates on the derivative by the same way because we can derive $A$ and we obtain the same type of quantities to estimate.
\end{prooff}

Now  we compute the boundary symbol obtained in \eqref{eq traces 0} and \eqref{eq bord 2}.
\begin{lemma}\label{lem sym 1}
 Let $k=0,1$ and  $\Im \rho_1>0$, $\Im \rho_2<0$, we have 
$$\left(  \int_\R e^{i\frac{x_n}{h}\xi_n} \frac{\xi_n^k}{(\xi_n-\rho_1)(\xi_n-\rho_2)}d\xi_n\right) _{|x_n=0}=2i\pi \frac{\rho_1^k}{\rho_1-\rho_2},
$$
\end{lemma}

\begin{prooff}
As in the proof of Lemma~\ref{lem trace}, we can integrate on $\Gamma$ 
and this integral is equal to
$2i\pi$ times the residu at $\rho_1$. It is easy to see that the residu is $\frac{\rho_1^k}{\rho_1-\rho_2}$.
\end{prooff}

\begin{lemma}\label{lem sym 2}
Let $k=0,1$ and $\Im \lambda_1>0$, $\Im \lambda_2<0$, $\Im \rho_1>0$, $\Im \rho_2<0$, we have 
$$\left(  \int_\R e^{i\frac{x_n}{h}\xi_n} \frac{\xi_n^k}{(\xi_n-\lambda_1)(\xi_n-\lambda_2)(\xi_n-\rho_1)(\xi_n-\rho_2)}d\xi_n\right) _{|x_n=0}=2i\pi A_k ,
$$
where 
\begin{equation*}
A_k=\left\lbrace 
	\begin{array}{l}
		\dfrac{\rho_2-\rho_1+\lambda_2-\lambda_1}{(\lambda_1-\lambda_2)(\lambda_1-\rho_2)(\rho_1-\lambda_2)(\rho_1-\rho_2)} \text{ if } k=0 \\[10 pt]
		\dfrac{\lambda_2\rho_2-\lambda_1\rho_1}{(\lambda_1-\lambda_2)(\lambda_1-\rho_2)(\rho_1-\lambda_2)(\rho_1-\rho_2)} \text{ if } k=1.
	\end{array}
\right. 
\end{equation*}

\end{lemma}

\begin{prooff}
 Clearly  both sides of the equality are continuous with respect $(\lambda_1,\rho_1)$ then it is sufficient to prove the case $\lambda_1\not=\rho_1 $.

As in the proof of Lemma~\ref{lem trace}, we can integrate on $\Gamma$ and the result is
$2i\pi$ times the sum of the residues in half plane $\Im z>0$. We obtain
\begin{align}
 &\frac{\lambda_1^k}{(\lambda_1-\lambda_2)(\lambda_1-\rho_1)(\lambda_1-\rho_2)}+\frac{\rho_1^k}{(\rho_1-\lambda_1)(\rho_1-\lambda_2)(\rho_1-\rho_2)}\notag\\
&\qquad\qquad\qquad\qquad\qquad\qquad\qquad\qquad =\frac{\lambda_1^k(\rho_1-\lambda_2)(\rho_1-\rho_2)-\rho_1^k(\lambda_1-\lambda_2)(\lambda_1-\rho_2)}{(\lambda_1-\lambda_2)(\lambda_1-\rho_1)(\lambda_1-\rho_2)(\rho_1-\lambda_2)(\rho_1-\rho_2)}.\notag
\end{align}
Clearly the numerator is null if $\rho_1=\lambda_1$. By a straightforward computation, if $k=0$ the numerator is $(\lambda_1-\rho_1)(\rho_2-\rho_1+\lambda_2-\lambda_1)$.
If $k=1$ the numerator is $ (\lambda_1-\rho_1)(\lambda_2\rho_2-\lambda_1\rho_1) $. This gives the Lemma.
\end{prooff}


\begin{thebibliography}{}

\bibitem{Ag65} Agmon Shmuel, Lectures on elliptic boundary value problems.
     Van Nostrand Mathematical  Studies, No. 2, 1965.   

\bibitem{Bo71} Boutet de Monvel  Louis,
Boundary problems for pseudo-differential operators. Acta Math. \textbf{126} (1971) 11--51.


\bibitem{CCC08} Cakoni Fioralba, {\c{C}}ay{\"o}ren Mehmet,  Colton
              David, Transmission eigenvalues and the nondestructive testing of  dielectrics. Inverse Problems  \textbf{24} (2008) 065016.

\bibitem{CCH10} Cakoni Fioralba, Colton David, Haddar Houssem, The interior transmission problem for regions with cavities. SIAM J. Math. Anal. \textbf{42} (2010) 145--162.


\bibitem{CGH10} Cakoni Fioralba, Gintides Drossos, Haddar Houssem, The existence of an infinite discrete set of transmission eigenvalues. SIAM J. Math. Anal. \textbf{42}  (2010) 237--255.

\bibitem{CH12} Cakoni Fioralba, Haddar Houssem, Transmission Eigenvalues in Inverse Scattering Theory. http://hal.inria.fr/hal-00741615

\bibitem{CKP89} Colton David, Kirsch  Andreas, P{\"a}iv{\"a}rinta
              Lassi, Far-field patterns for acoustic waves in an inhomogeneous medium. SIAM J. Math. Anal. \textbf{20} (1989) 1472--1483.



\bibitem{CK92} Colton  David, Kress  Rainer, Inverse acoustic and electromagnetic scattering theory. Applied Mathematical Sciences, \textbf{93}, {Springer-Verlag}, {Berlin}, {1992}.  


\bibitem{CPS07}Colton  David, P{\"a}iv{\"a}rinta  Lassi, Sylvester
              John, The interior transmission problem. Inverse Probl. Imaging \textbf{1} (2007) 13--28.

\bibitem{CR12} Cornilleau Pierre, Robbiano Luc, Carleman estimates for the Zaremba Boundary Condition and Stabilization of Waves. http://hal.archives-ouvertes.fr/hal-00634867/fr/

\bibitem{Da99} Davies E. B., Semi-classical states for non-self-adjoint {S}chr\"odinger
              operators. Communications in Mathematical Physics \textbf{200} (1999) 35--41.

\bibitem{DSZ04} Dencker  Nils,  Sj{\"o}strand  Johannes, Zworski  Maciej,
     Pseudospectra of semiclassical (pseudo-) differential
              operators. Comm. Pure Appl. Math. \textbf{57} {(2004)} {384--415.}

\bibitem{HKOP10} Hitrik Michael, Krupchyk  Katsiaryna, Ola Petri,
              P{\"a}iv{\"a}rinta Lassi. Transmission eigenvalues for operators with constant
              coefficients. SIAM J. Math. Anal. \textbf{42} (2010) 2965--2986.

\bibitem{HKOP11} Hitrik Michael, Krupchyk Katsiaryna,  Ola Petri,
              P{\"a}iv{\"a}rinta Lassi, Transmission eigenvalues for elliptic operators. SIAM J. Math. Anal. \textbf{43} (2011) 2630--2639.

\bibitem{Hobo83} H\"{o}rmander Lars, The analysis of linear partial differential operators, \textit{Vol.I-IV}. Grundlehren der Mathematischen Wissenschaften, \textbf{255-256}, \textbf{274-275}. Springer-Verlag, Berlin.  1983, 1985.

\bibitem{LV12} Lakshtanov E., Vainberg, B., Ellipticity in the {I}nterior {T}ransmission {P}roblem in {A}nisotropic {M}edia. SIAM J. Math. Anal. \textbf{44} (2012) 1165--1174.

\bibitem{LV11-arxiv} Lakshtanov E., Vainberg B. Remarks on interior transmission eigenvalues, Weyl formula and branching billiards. arxiv: arXiv:1112.0891

\bibitem{LV12-arxiv} Lakshtanov E., Vainberg B. Bounds on positive interior transmission eigenvalues. arxiv: 1206.3782v2.

\bibitem{LV12bis-arxiv} Lakshtanov E., Vainberg B. Applications of elliptic operator theory to the isotropic interior transmission eigenvalue problem. arxiv: arXiv:1212.6785

\bibitem{LR95} Lebeau Gilles. Robbiano Luc, Contr\^{o}le exact de l'\'equation de la chaleur. \textit{Comm. Partial Differential Equations}, \textbf{20} ( 1995) 335-356.

\bibitem{LR97} Lebeau Gilles, Robbiano Luc,  Stabilisation de l'\'equation des ondes par le bord. \textit{Duke Math. J.}, \textbf{86}( 1997) 465-491.

\bibitem{PS08} P{\"a}iv{\"a}rinta Lassi, Sylvester  John, Transmission eigenvalues. SIAM J. Math. Anal. \textbf{40} (2008) 738--753.

\bibitem{Sy12} Sylvester john, Discreteness of Transmission Eigenvalues via Upper Triangular Compact Operators. SIAM J.  Math.  Anal. \textbf{44} (2012) 341--354.

\bibitem{Zw01} Zworski Maciej, A remark on a paper of {E}. {B} {D}avies: ``{S}emi-classical states for non-self-adjoint {S}chr\"odinger operators''
              [{C}omm. {M}ath. {P}hys. {\bf 200} (1999), no. 1, 35--41;
              {MR}1671904 (99m:34197)]. Proceedings of the American Mathematical Society  \textbf{129} (2001) 2955--2957.
\end{thebibliography}
\end{document}